\documentclass[11pt]{article}
\usepackage{mathrsfs}
\usepackage{amssymb}
\usepackage{amsmath}
\usepackage{amsfonts}
\usepackage{graphicx}

\newtheorem{theorem}{Theorem}[section]
\newtheorem{remark}{Remark}[section]
\newtheorem{corollary}{Corollary}[section]

\date{}
\begin{document}

\title{Differential games of partial information forward-backward doubly stochastic differential equations
and applications
 \thanks{ This research project was
funded by PolyU research accounts 1-ZV1X and G-YH86 of Hong Kong,
the National Nature Science Foundation of China (11001156, 11071144,
11026125), the Nature Science Foundation of Shandong Province
(ZR2009AQ017), and Independent Innovation Foundation of Shandong
University (IIFSDU), China. }
 }
\author{Eddie C.M. HUI\thanks{Department of Building and Real Estate,
The Hong Kong Polytechnic University, Hung Hom, Kowloon, Hong Kong.
E-mail address: bscmhui@polyu.edu.hk }\quad  and\;
   Hua XIAO \thanks{Corresponding author.
School of Mathematics and Statistics, Shandong University at Weihai,
       Weihai 264209, China;
       School of Mathematics, Shandong University,
       Jinan 250100, China.
      E-mail address:  xiao$\_$\ hua@sdu.edu.cn}
 }

 \maketitle
 \noindent{\bf Abstract}\; This paper is concerned with a new type of differential game problems of forward-backward stochastic systems. There are three distinguishing features: Firstly, our game systems are forward-backward doubly stochastic differential equations, which is a class of more general game systems than other forward-backward stochastic game systems without doubly stochastic terms; Secondly, forward equations are directly related to backward equations at initial time, not terminal time; Thirdly, the admissible control is required to be adapted to a sub-information of the full information generated by the underlying Brownian motions. We give a necessary  and a sufficient conditions for both an equilibrium point of nonzero-sum games and a saddle point of zero-sum games.
Finally, we work out an example of linear-quadratic nonzero-sum differential games to
illustrate the theoretical applications. Applying some stochastic filtering techniques, we obtain the explicit expression of the equilibrium point.

\vspace{0.3cm}

\vspace{0.3cm}

\section{Introduction}

Game theory is a useful tool which helps us understand economic, social, political,
and biological phenomena. Stochastic differential game problems also attract more and more research attentions, and are used widely in other social and behavioral sciences. When we study stochastic differential games of backward doubly stochastic differential equations (BDSDEs, for short), doubly stochastic Hamiltonian systems with boundary conditions appear naturally whose dynamics are described by initial coupled forward-backward doubly stochastic differential equations (FBDSDEs, for short). To illustrate this, we introduce an example of linear quadratic (LQ, for short) nonzero-sum differential games of BDSDEs with partial information which  motivates us to initiate a study of stochastic differential games of initial coupled FBDSDEs with partial information. We now explain this in more
detail.

Let $T$ be a fixed constant and $\Big(\Omega, \mathscr{F},
 P\Big)$ be a complete
filtered probability space, on which two mutually independent
standard Brownian motions $B(\cdot)\in\mathrm{R}^{l}$ and
$W(\cdot)\in\mathrm{R}^{d}$ are defined. Let $\mathcal{N}$ denote
the class of $P$-null sets of $\mathcal{F}$. For each $t\in[0, T],$
we define
$$\mathcal{F}_{t}\doteq \mathcal{F}^{W}_{t}\vee \mathcal{F}^{B}_{t, T},$$
where $\mathcal{F}^{W}_{t}=\mathcal{N}\vee \sigma\{W(r)-W(0): 0\leq
r\leq t\}$ and $\mathcal{F}^{B}_{t, T}=\mathcal{N}\vee
\sigma\{B(r)-B(t): t\leq r\leq T\}.$ Note that the set
$\mathcal{F}_{t}, t\in[0, T]$ is neither increasing nor decreasing,
so it does not constitute a filtration. We denote by
$\mathcal{L}^{p}_{T}(\Omega; \mathcal{S})$ all class of
$\mathcal{F}_{T}$-measurable random variables $\{\xi:
\Omega\longrightarrow \mathcal{S}\}$ satisfying $E|\xi|^{p}<\infty$,
by $\mathcal{L}^{p}_{\mathcal{F}_{t}}(0, T; \mathcal{S})$ all class
of $\mathcal{F}_{t}$-adapted stochastic processes $\{x(t): [0,
T]\times \Omega\longrightarrow \mathcal{S}\}$ satisfying
$\mathbb{E}\big[\int_{0}^{T}|x(t)|^{p}dt\big]<+\infty$. If there is
no risk of confusion, we write
$\mathcal{L}^{p}_{T}=\mathcal{L}^{p}_{T}(\Omega;\mathcal{S}),$
$\mathcal{L}^{p}_{\mathcal{F}_{t}}=\mathcal{L}^{p}_{\mathcal{F}_{t}}(0,
T; \mathcal{S})$. The processes $v_{1}(t)=v_{1}(t, \omega)$ and
$v_{2}(t)=v_{2}(t, \omega)$ are our open-loop control processes. Let $U_{i}$
be a nonempty convex subset of $\mathrm{R}^{k_{i}}$ ( $i=1, 2$). In many cases in which the full information $\mathcal{F}_{t}$ is inaccessible
for players, ones can only observe a partial information. For this,
we denote the set of all
open-loop admissible controls for the player $i$ by
\begin{align*}
  \mathcal{U}_{i}=\Big\{v_{i}(\cdot): [0, T]\times \Omega\longrightarrow U_{i}\big| v_{i}(\cdot)
   \hbox{ is } \mathcal{E}_{t}\hbox{-adapted and satisfies}\; \mathbb{E}\int_{0}^{T}|v_{i}(t)|^{2}dt<\infty\Big\},
\end{align*}
where $i=1, 2,$
$\mathcal{E}_{t}$ is an available sub-information of full information $\mathcal{F}_{t}$ for players, i.e.
$$\mathcal{E}_{t}\subseteq \mathcal{F}_{t},\quad\quad \hbox{for all}\ t.$$
For example, $\mathcal{E}_{t}$ could be the $\delta$-delayed information defined by
$$\mathcal{E}_{t}=\mathcal{F}_{(t-\delta)^{+}},$$ where $\delta$ is a given positive constant delay.
Each element of $\mathcal{U}_{i}$ is called an admissible control
for Player $i$ on [0, T] $ (i=1, 2)$.
$\mathcal{U}_{1}\times\mathcal{U}_{2}$ is called the set of
open-loop admissible controls for the players.

 We consider the following 1-dimensional linear BDSDE
\begin{equation}
    \left\{
    \begin{aligned}\label{DG26}
    -dY(t)=&[A_{1}Y(t)+B_{1}Z(t)+C_{1}v_{1}(t)+D_{1}v_{2}(t)]dt\\
           &+[A_{2}Y(t)+B_{2}Z(t)+C_{2}v_{1}(t)+D_{2}v_{2}(t)]\hat{d}B(t)-Z(t)dW(t),\\
    Y(T)=& \xi,
    \end{aligned}
    \right.
\end{equation}
and the performance criterion, for $i=1, 2,$
\begin{equation}\label{DG27}
\begin{aligned}
    J_{i}(v_{1}(\cdot), v_{2}(\cdot))=&-\frac{1}{2}\mathbb{E}\Bigg\{\langle F_{i1}Y(0), Y(0)\rangle+\int_{0}^{T}\Big[\langle F_{i2}Y(t), Y(t)\rangle\\
                          &\hspace{1cm}+\langle  F_{i3}Z(t), Z(t)\rangle+\langle F_{i4}v_{1}(t), v_{1}(t)\rangle+\langle F_{i5}v_{2}(t), v_{2}(t)\rangle\Big]dt\Bigg\},
\end{aligned}
\end{equation}
where the integral with respect to $\hat{d}B(t)$ is a
"backward It\^{o} integral" and the integral with respect to
$dW(t)$ is a standard forward It\^{o} integral. These are two types
of particular cases of the It\^{o}-Skorohod integral (see Nualart
and Pardoux\cite{NP1988}). The extra noise $\{B(t)\}$ can be
considered some extra information that can not be detected in
practice, such as in a derivative security market, but is valuable
to the partial investors.

If we let $A_{2}, B_{2}, C_{2}, D_{2}\equiv 0$, then equation \eqref{DG26} is reduced to a general backward stochastic differential equation (BSDE, for short) of Pardoux-Peng's type (see Pardoux and Peng\cite{PP1990}).
For simplicity, we assume temporarily that $\xi\in\mathcal{L}^{2}_{T}(\Omega,\mathrm{R}^{1}),$ all coefficients in \eqref{DG26} and \eqref{DG27} are 1-dimensional, $l=d=1,$
 $F_{i1}, F_{i2},  F_{i3}\geq 0$, $F_{i4}, F_{i5}>0$.

Our aim is to seek an equilibrium point $(u_{1}(\cdot), u_{2}(\cdot))\in
\mathcal{U}_{1}\times\mathcal{U}_{2}$,
for all  $(v_{1}(\cdot), v_{2}(\cdot))\in
\mathcal{U}_{1}\times\mathcal{U}_{2}$, such that
\begin{equation*}
\left\{
\begin{aligned}
    J(u_{1}(\cdot), u_{2}(\cdot))&\geq J(v_{1}(\cdot), u_{2}(\cdot)),\\
     J(u_{1}(\cdot), u_{2}(\cdot))&\geq J(u_{1}(\cdot), v_{2}(\cdot)).\\
\end{aligned}
\right.
\end{equation*}
We call it an LQ nonzero-sum differential game of BDSDE and denote it by {\it Problem (LQNZB)}.

Applying Theorem 4.1 in Han et al.\cite{HPW2010}, we conclude that the equilibrium point must satisfy the following form:
\begin{equation}\label{DG35}
\left\{
\begin{aligned}
  &u_{1}(t)=-\mathbb{E}\big[F_{14}^{-1}\big(C_{1}y_{1}(t)+C_{2}z_{1}(t)\big)\big|\mathcal{E}_{t}\big],\\
  &u_{2}(t)= \mathbb{E}\big[F_{25}^{-1}\big(D_{1}y_{2}(t)+D_{2}z_{2}(t)\big)\big|\mathcal{E}_{t}\big],
\end{aligned}
\right.
\end{equation}
where $(y_{i}, z_{i})\ (i=1,2)$ is the solution of the following initial coupled FBDSDE:
\begin{equation}\label{DG28}
    \left\{
    \begin{aligned}
    -dY(t)=&\Big[A_{1}Y(t)+B_{1}Z(t)-C_{1}F_{14}^{-1}\big(C_{1}y_{1}(t)+C_{2}z_{1}(t)\big)\\
           & -D_{1}F_{25}^{-1}\big(D_{1}y_{2}(t)+D_{2}z_{2}(t)\big)\Big]dt\\
           &+\Big[A_{2}Y(t)+B_{2}Z(t)-C_{2}F_{14}^{-1}\big(C_{1}y_{1}(t)+C_{2}z_{1}(t)\big)\\
           & -D_{2}F_{25}^{-1}\big(D_{1}y_{2}(t)+D_{2}z_{2}(t)\big)\hat{d}B(t)-Z(t)dW(t),\\
    dy_{i}(t)=& \big(A_{1}y_{i}(t)+A_{2}z_{i}(t)+F_{i2}Y(t)\big)dt+\big(B_{1}y_{i}(t)+B_{2}z_{i}(t)+ F_{i3}Z(t)\big)dW(t)\\
          & -z_{i}(t)\hat{d}B(t),\\
    Y(T)=& \xi,\; y_{i}(0)=F_{i1}Y(0).
    \end{aligned}
    \right.
\end{equation}

We see that the equation for $Y(\cdot)$ is backward
(since it is given the final datum which is an $\mathcal{F}_{T}$-measurable random variable), the equation for $y_{i}(\cdot)$ is forward (since it is given the initial datum which is directly related with the backward solution $Y(\cdot)$ at initial time). Further, the backward equation is ``forward" with respect to the backward stochastic integral $\hat{d}B(t),$ as well as ``backward" with respect to the forward stochastic  integral $dW(t);$ the coupled forward equation is ``backward" with respect to the backward stochastic  integral $\hat{d}B(t),$ as well as ``forward" with respect to the forward stochastic  integral $dW(t).$  Equation \eqref{DG28} is exactly the type of time-symmetric forward-backward stochastic differential equations (FBSDEs, for short) introduced by Peng and Shi\cite{PS2003}. There is a small difference between equation \eqref{DG28} and FBSDE in Peng and Shi\cite{PS2003}: the former is initial coupled, but the latter is terminal coupled. So we call equation \eqref{DG28} an initial coupled linear FBDSDE.  In addition, the candidate equilibrium point denoted by \eqref{DG35} involves the available sub-information $\mathcal{E}_{t}$ of full information $\mathcal{F}_{t}$ for players. Thus, a type of initial coupled  FBDSDE naturally appears when we study {\it Problem (LQNZB)}. Since this type of FBDSDEs possess fine dynamics and can be reduced to FBSDEs or BDSDEs or BSDEs, one could not help thinking about differential game problems for initial coupled FBDSDEs under partial information.

Pardoux\cite{Pardoux1979} generalized the classical Feynman-Kac formula and provided a probabilistic representation for solutions of linear parabolic stochastic partial differential equations (SPDEs, for short). By introducing originally BDSDEs, which is a new class of BSDEs and covers the results of Pardoux and Peng\cite{PP1990}, Pardoux and Peng\cite{PP1994} produced a probabilistic representation of certain quasi-linear SPDEs as an extension to the Feynman-Kac formula for linear SPDEs. In general, that a forward SDE of It\^{o}'s type couples a backward SDE of Pardoux-Peng's type, which maybe couple each other at initial conditions or terminal conditions, constitutes an initial or terminal coupled FBSDE. The theory of FBSDEs has received considerable research attention in recent years. For more information on the solvability of FBSDEs and corresponding optimal control problems  with full or partial information, see e.g. Antonelli\cite{Antonelli1993}, Hu and Peng\cite{HP1995}, Ma, Protter and Yong\cite{MPY1994}, Meng\cite{Meng2009},
{\O}ksendal and Sulem\cite{OS2010}, Peng and Shi\cite{PS2000}, Peng and Wu\cite{PW1999}, Shi and Wu\cite{SW2006,SW-appear}, Wang and Wu\cite{WW2008, WW2009},  Wu\cite{Wu2005,Wu2010},  specially the monographs by Ma and Yong\cite{MY1999} and Yong and Zhou\cite{YZ1999}, etc.

There is a few literature on differential games of BSDEs and FBSDEs. Yu and Ji\cite{YJ2008} obtained an existence and uniqueness result for an initial coupled FBSDE under some monotone conditions, applied it to backward linear-quadratic nonzero-sum stochastic differential game problem and got the explicit form of a Nash equilibrium point. Wang and Yu\cite{WY2010} established a necessary and a sufficient conditions for an equilibrium point of nonzero-sum differential game of BSDEs and applied them to study a financial problem. Zhang\cite{Zhang2009} extended the result of Yu and Ji\cite{YJ2008} to the case where BSDEs are driven by both Brownian motion and Poisson random measure. Wang and Yu\cite{WY2011} recently generalize the results of \cite{WY2010} to partial information differential games and obtain the corresponding maximum principle and verification theorem, and they also apply the theoretical results to study LQ differential games and financial problem. Yu\cite{Yu2012} mainly studied the LQ optimal control and nonzero-sum differential game of FBSDE. Hui and Xiao\cite{HX2011} investigated differential games of FBSDEs, and established the maximum principle and verification theorem for both an equilibrium point of nonzero-sum cases and a saddle point of zero-sum cases. Meng\cite{Meng2010} discussed the partial information zero-sum differential games of fully coupled FBSDEs.

Han et al.\cite{HPW2010} investigated the optimal control for BDSDEs and obtain a stochastic maximum principle of the optimal control. In \cite{PS2003}, Peng and Shi established the existence and uniqueness results of terminal coupled FBDSDEs under certain monotonicity assumptions. Zhu et al.\cite{ZSG2009} relaxed the monotonicity assumptions and allowed the case of different dimensions between forward equations and backward equations, compared with the results in Peng and Shi\cite{PS2003}.
Zhang and Shi\cite{ZS2010} studied the optimal control of fully terminal-coupled FBDSDEs and obtained the maximum principle in the global form, where the control variables can enter into the diffusion coefficients and the control domain need not be convex. Using the solution of FBDSDEs, Zhang and Shi also got the explicit
form of open-loop Nash equilibrium point for nonzero-sum stochastic differential games of only forward doubly stochastic differential equations.

However, none of the works mentioned above deals with differential games of initial coupled FBDSDEs with full information or partial information.
In Section 2, we formulate the
zero-sum and nonzero-sum games of initial-coupled FBDSDEs with partial information. In Section 3,  we are devoted to proving a
maximum principle and a verification theorem for both an equilibrium point of nonzero-sum games and a saddle point of zero-sum
games.  In Section 4, an
example of a nonzero-sum differential game is worked out to
illustrate theoretical applications. In terms of maximum principle and verification theorem,
the explicit expression of an equilibrium point is obtained. Finally,
we give some concluding remarks.

\section{Formulation of the problem}

We introduce the mappings
\begin{equation*}
    \begin{aligned}
&f: [0,T]\times\mathrm{R}^{n}\times\mathrm{R}^{n\times
l}\times\mathrm{R}^{m}\times\mathrm{R}^{m\times d}\times
    U_{1}\times U_{2}\rightarrow\mathrm{R}^{n},\\
&\bar{f}: [0,T]\times\mathrm{R}^{n}\times\mathrm{R}^{n\times
l}\times\mathrm{R}^{m}\times\mathrm{R}^{m\times d}\times
    U_{1}\times U_{2}\rightarrow\mathrm{R}^{n \times d},\\
&g:[0,T]\times\mathrm{R}^{m}\times\mathrm{R}^{m\times
d}\times U_{1}\times U_{2}\rightarrow\mathrm{R}^{m},\\
&\bar{g}:[0,T]\times\mathrm{R}^{m}\times\mathrm{R}^{m\times
d}\times U_{1}\times U_{2}\rightarrow\mathrm{R}^{m\times l},\\
&\phi: \mathrm{R}^{m}\rightarrow\mathrm{R}^n,\quad\quad \varphi,
\varphi_{i}: \mathrm{R}^{m}\rightarrow\mathrm{R}^1,
    \quad\quad \gamma, \gamma_{i}: \mathrm{R}^{n}\rightarrow\mathrm{R}^1,
    \quad \\
&l, l_{i}: [0,T]\times\mathrm{R}^{n}\times\mathrm{R}^{n\times
l}\times\mathrm{R}^{m}\times\mathrm{R}^{m\times d}\times U_{1}\times
U_{2}\rightarrow\mathrm{R}^{1}\ (i=1, 2).
    \end{aligned}
\end{equation*}

\noindent Assumption (H1): For any $(y, z, Y, Z, v_{1},
v_{2})\in\mathrm{R}^{n}\times\mathrm{R}^{n\times
l}\times\mathrm{R}^{m}\times\mathrm{R}^{m\times d}\times
    U_{1}\times U_{2},$ we assume that
\begin{align*}
    &f(\cdot, y, z, Y, Z, v_{1}, v_{2})\in\mathcal{L}^{2}_{\mathcal{F}_{t}}(0, T;
 \mathcal{R}^{n}),\quad
 \bar{f}(\cdot, y, z, Y, Z, v_{1}, v_{2})\in\mathcal{L}^{2}_{\mathcal{F}_{t}}(0, T;
 \mathcal{R}^{n\times d}),\\
&g(\cdot,  Y, Z, v_{1},
v_{2})\in\mathcal{L}^{2}_{\mathcal{F}_{t}}(0, T;
 \mathcal{R}^{m}),\quad\quad\quad
 \bar{g}(\cdot,  Y, Z, v_{1},
v_{2})\in\mathcal{L}^{2}_{\mathcal{F}_{t}}(0, T;
 \mathcal{R}^{m\times l}).
\end{align*}
We assume moreover that $f, \bar{f}, g$ and $\bar{g}$   are
continuously differentiable with respect to $ (y, z, Y, Z,$ $v_{1},
v_{2})$ and their derivatives with respect to $ (y, z, Y, Z, v_{1},
v_{2})$ are continuous and uniformly bounded. $l, l_{1}, l_{2},
\varphi, \varphi_{1},  \varphi_{2}, \gamma, \gamma_{1}$ and $
\gamma_{2}$ are continuously differential with respect to $ (y, z,
Y, Z, v_{1},$ $v_{2})$ and their derivatives with respect to $ (y, z,
Y, Z, v_{1}, v_{2})$ are continuous and bounded by
$K(1+|y|+|z|+|Y|+|Z|+|v_{1}|+|v_{2}|)$. There exists constants $k>0$
and $0<c<1$ such that 
\begin{align*}&|\bar{f}(t, y_{1}, z_{1}, Y_{1}, Z_{1}, u_{1}, u_{2})-\bar{f}(t, y_{2}, z_{2},
Y_{2}, Z_{2}, v_{1}, v_{2})|^{2}\\
\leq& k(|y_{1}-y_{2}|^{2}+|Y_{1}-Y_{2}|^{2}
+|Z_{1}-Z_{2}|^{2}+|u_{1}-v_{1}|^{2}+|u_{2}-v_{2}|^{2})+c|z_{1}-z_{2}|^{2},\\
&|\bar{g}(t, Y_{1}, Z_{1}, u_{1}, u_{2})-\bar{g}(t, Y_{2}, Z_{2},
v_{1}, v_{2})|^{2}\\
\leq&k(|Y_{1}-Y_{2}|^{2}+|u_{1}-v_{1}|^{2}+|u_{2}-v_{2}|^{2})+c|Z_{1}-Z_{2}|^{2}),
\end{align*}
for all $ (y_{1}, z_{1}, Y_{1}, Z_{1}, u_{1}, u_{2}),$ $ (y_{2},
z_{2}, Y_{2}, Z_{2}, v_{1},
v_{2})\in\mathrm{R}^{n}\times\mathrm{R}^{n\times
l}\times\mathrm{R}^{m}\times\mathrm{R}^{m\times d}\times
    U_{1}\times U_{2}.$

\par
In the following, we specify the problems of nonzero-sum
 and zero-sum differential games of
forward-backward doubly stochastic systems, respectively. For
simplicity, we denote them by {\it Problem (NZSG)} and {\it Problem
(ZSG)}, respectively.
\par
Consider an FBDSDE
 {\setlength\arraycolsep{1pt}
\begin{equation}\label{DG14}
\left\{\begin{aligned}
-dY^{v_{1},v_{2}}(t)=&\ g(t,  Y^{v_{1},v_{2}}(t), Z^{v_{1},v_{2}}(t), v_{1}(t), v_{2}(t))dt\\
&+\bar{g}(t,  Y^{v_{1},v_{2}}(t), Z^{v_{1},v_{2}}(t), v_{1}(t), v_{2}(t))\hat{d}B(t)-Z^{v_{1},v_{2}}(t)dW(t),\\
dy^{v_{1},v_{2}}(t)=&\ f(t, y^{v_{1}, v_{2}}(t),  z^{v_{1}, v_{2}}(t), Y^{v_{1},v_{2}}(t), Z^{v_{1},v_{2}}(t), v_{1}(t), v_{2}(t))dt\\
                &+\bar{f}(t, y^{v_{1},v_{2}}(t), z^{v_{1}, v_{2}}(t), Y^{v_{1},v_{2}}(t), Z^{v_{1},v_{2}}(t), v_{1}(t), v_{2}(t))dW(t)\\
                &-z^{v_{1},v_{2}}(t)\hat{d}B(t),\\
Y^{v_{1},v_{2}}(T)=&\ \xi,\quad \quad y^{v_{1},v_{2}}(0)=
\phi(Y^{v_{1},v_{2}}(0)),\quad \quad 0\leq t\leq
             T.
\end{aligned}
\right.
\end{equation}
}
Under the assumption (H1), there exists a
unique solution $\big(y^{v_{1}, v_{2}}(\cdot), z^{v_{1},
v_{2}}(\cdot), Y^{v_{1}, v_{2}}(\cdot), Z^{v_{1},
v_{2}}(\cdot)\big)\in\mathcal{L}^{2}_{\mathcal{F}_{t}}(0,T;
\mathrm{R}^{n})\times\mathcal{L}^{2}_{\mathcal{F}_{t}}(0,T;
\mathrm{R}^{n\times l})\times\mathcal{L}^{2}_{\mathcal{F}_{t}}(0,T;
\mathrm{R}^{m})\times\mathcal{L}^{2}_{\mathcal{F}_{t}}(0,T;
\mathrm{R}^{m\times d})$ to equation \eqref{DG14} for any $
(v_{1}(\cdot), v_{2}(\cdot))$ $\in\mathcal{U}_{1}\times\mathcal{U}_{2}$
(see Pardoux and Peng\cite{PP1994}). In the case where equation \eqref{DG14} does not involve the term of backward It\^{o}'s integral, i.e. $\bar{g}\equiv 0$, $f$ and $\bar{f}$
are independent of $z^{v_{1}, v_{2}}$, the game systems will be reduced to the initial coupled FBSDEs which has been studied by Xiao and Wang\cite{XWJAMC}.
In the case where equation \eqref{DG14} does not involve the forward equation, i.e. $f=\bar{f}=\phi\equiv 0$, the game systems will be reduced to the BDSDEs which has been investigated by Han et al.\cite{HPW2010}.
In the case where equation \eqref{DG14} does not involve both the term of backward It\^{o}'s integral and the forward equation, the game systems will be reduced to the BSDEs which has been investigated by Wang and Yu\cite{WY2010,WY2011} and Yu and Ji\cite{YJ2008}.

\par
Consider a performance criterion
\begin{align}\label{DG13}
J_{i}(v_{1}(\cdot), v_{2}(\cdot))=\mathbb{E}\Big[&\int_{0}^{T}l_{i}
\big(t, y^{v_{1},v_{2}}(t),  z^{v_{1},v_{2}}(t), Y^{v_{1},v_{2}}(t),
Z^{v_{1},v_{2}}(t), v_{1}(t),
v_{2}(t)\big)dt\nonumber\\
&+\varphi_{i}(Y^{v_{1},v_{2}}(0))\Big]+\gamma_{i}(y^{v_{1},v_{2}}(T))
\end{align}
with $l_{i}(\cdot, y^{v_{1},v_{2}}(\cdot), z^{v_{1},v_{2}}(\cdot),
Y^{v_{1},v_{2}}(\cdot), Z^{v_{1},v_{2}}(\cdot), v_{1}(\cdot),
v_{2}(\cdot))\in\mathcal{L}^{1}_{\mathcal{F}_{t}}(0,T; \mathrm{R}) $
and $\varphi_{i}\in\mathcal{L}^{1}$ $(0, T; \mathrm{R})$ for any $
(v_{1}(\cdot), v_{2}(\cdot))\in\mathcal{U}_{1}\times\mathcal{U}_{2}
\ (i=1, 2)$, and
\begin{align}\label{DG12}
J(v_{1}(\cdot), v_{2}(\cdot))=\mathbb{E}\Big[\int_{0}^{T}&l \big(t,
y^{v_{1},v_{2}}(t),  z^{v_{1},v_{2}}(t), Y^{v_{1},v_{2}}(t),
Z^{v_{1},v_{2}}(t), v_{1}(t),
v_{2}(t)\big)dt\nonumber\\
&+\varphi(y^{v_{1},v_{2}}(T))\Big]+\gamma(Y^{v_{1},v_{2}}(0))
\end{align}
with $l(\cdot, y^{v_{1},v_{2}}(\cdot), Y^{v_{1},v_{2}}(\cdot),
Z^{v_{1},v_{2}}(\cdot), v_{1}(\cdot),
v_{2}(\cdot))\in\mathcal{L}^{1}_{\mathcal{F}_{t}}(0,T; \mathrm{R}) $
and $\varphi\in\mathcal{L}^{1}$ $(0, T; \mathrm{R})$ for any $
(v_{1}(\cdot),
v_{2}(\cdot))\in\mathcal{U}_{1}\times\mathcal{U}_{2}$. We note that
\eqref{DG13} and \eqref{DG12} are well posed. There are two players
$i_{1}$ and $i_{2}$. Player $i_{1}$ controls $v_{1}$ and Player
$i_{2}$ controls $v_{2}$.

{\it Problem (NZSG):}\; Find $(u_{1}(\cdot),
u_{2}(\cdot))\in\mathcal{U}_{1}\times\mathcal{U}_{2} $ such that
\begin{equation}\label{DG11}
\left\{
\begin{aligned}
 J_{1}(u_{1}(\cdot), u_{2}(\cdot))\geq J_{1}(v_{1}(\cdot), u_{2}(\cdot)),\\
J_{2}(u_{1}(\cdot), u_{2}(\cdot))\geq J_{2}(u_{1}(\cdot),
v_{2}(\cdot)),
\end{aligned} \right.
\end{equation}
for all $(v_{1}(\cdot),
v_{2}(\cdot))\in\mathcal{U}_{1}\times\mathcal{U}_{2} $.
We call $(u_{1}(\cdot), u_{2}(\cdot))$ an
open-loop equilibrium point of {\it Problem (NZSG)} (if it does exist). It is
easy to see that the existence of an open-loop equilibrium point
implies
\begin{equation*}
\left\{
\begin{aligned}
 J_{1}(u_{1}(\cdot), u_{2}(\cdot))=\sup_{v_{1}(\cdot)\in\mathcal{U}_{1}} J_{1}(v_{1}(\cdot), u_{2}(\cdot)),\\
J_{2}(u_{1}(\cdot),
u_{2}(\cdot))=\sup_{v_{2}(\cdot)\in\mathcal{U}_{2}}
J_{2}(u_{1}(\cdot), v_{2}(\cdot)).
\end{aligned} \right.
\end{equation*}

{\it Problem (ZSG):}\; Find
$(u_{1}(\cdot), u_{2}(\cdot))\in\mathcal{U}_{1}\times\mathcal{U}_{2}
$ such that
\begin{equation}\label{DG10}
\begin{aligned}
 J(u_{1}(\cdot), v_{2}(\cdot))\leq J(u_{1}(\cdot), u_{2}(\cdot))\leq J(v_{1}(\cdot),
 u_{2}(\cdot)),
\end{aligned}
\end{equation}
for all $(v_{1}(\cdot),
v_{2}(\cdot))\in\mathcal{U}_{1}\times\mathcal{U}_{2} $.
We call $(u_{1}(\cdot), u_{2}(\cdot))$ an
open-loop saddle point of {\it Problem (ZSG)} (if it exists). In
fact the existence of an open-loop saddle point implies
\begin{align*}
 J(u_{1}(\cdot), u_{2}(\cdot))&=\sup_{v_{2}(\cdot)\in \mathcal{U}_{2}}\Big(\inf_{v_{1}(\cdot)\in \mathcal{U}_{1}}J\big(v_{1}(\cdot),
    v_{2}(\cdot)\big)\Big)\\
                           &=\inf_{v_{1}(\cdot)\in \mathcal{U}_{1}}\Big(\sup_{v_{2}(\cdot)\in
                           \mathcal{U}_{2}}J\big(v_{1}(\cdot),
    v_{2}(\cdot)\big)\Big).
\end{align*}
We shall verify this point in Theorem 3.5 (iii).

\section{Differential games of FBDSDEs}
\subsection{Nonzero-sum case}

Suppose $(u_{1}(\cdot), u_{2}(\cdot))$ is an equilibrium point of
{\it Problem (NZSG)} with the trajectory $\Big(y(\cdot),$ $z(\cdot),
Y(\cdot), Z(\cdot)\Big)$ of (\ref{DG14}). For all $ t\in [0, T],$ let
$v_{i}(t)\in U_{i}$ be such that
$u_{i}(\cdot)+v_{i}(\cdot)\in\mathcal{U}_{i}\ (i=1,2).$ Notice that
$\mathcal{U}_{i}$ is convex, then for $0\leq\epsilon, \rho\leq 1, i=1,2$,
$$u_{1\epsilon}(t)=u_{1}(t)+\epsilon v_{1}(t)\in\mathcal{U}_{1},\ u_{2\rho}(t)=u_{2}(t)+\rho v_{2}(t)\in\mathcal{U}_{2},  \quad0\leq
t\leq T.$$
For simplicity, we denote
\begin{align*}
&f(t)=f\big(t, y(t), z(t), Y(t), Z(t), u_{1}(t), u_{2}(t)\big),\\
&g(t)=g\big(t, Y(t), Z(t), u_{1}(t), u_{2}(t)\big),\\
&Y^{u_{1\epsilon }}(t)=Y^{(u_{1}+\epsilon v_{1},
u_{2})}(t),\;Y^{u_{2\rho }}(t)=Y^{(u_{1}, u_{2}+\rho v_{2})}(t),\\
&h_{i}(\epsilon, \rho)=J_{i}(u_{1}+\epsilon v_{1}, u_{2}+\rho
v_{2}),
\end{align*}
define the processes $$
\hat{Y}^{1}(t)=\frac{d}{d\epsilon}\;Y^{u_{1\epsilon
}}(t)|_{\epsilon=0}, \;\;\hat{Y}^{2}(t)=\frac{d}{d\rho}\;Y^{u_{2\rho
}}(t)|_{\rho=0},
$$
and make the similar notations for $\bar{f}, \bar{g}, l_{i},
\hat{y}^{i}, \hat{z}^{i},\hat{Z}^{i}, i=1,2$. For $i=1,2,$ we have
the following variational equations:
\begin{equation*}
\left\{
\begin{aligned}
-d\hat{Y}^{i}(t)=&\hat{g}^{i}(t)dt+\hat{\bar{g}}^{i}(t)\hat{d}B(t)-\hat{Z}^{i}(t)dW(t),\\
d\hat{y}^{i}(t)=&\hat{f}^{i}(t)dt+\hat{\bar{f}}^{i}(t)dW(t)-\hat{z}^{i}(t)\hat{d}B(t),\\
\hat{Y}^{i}(T)=&0,\;\hat{y}^{i}(0)=\phi_{Y}(Y(0))\hat{Y}^{i}(0)
\end{aligned}
\right.
\end{equation*}
where
\begin{align*}
\hat{g}^{i}(t)=&g_{Y}(t)\hat{Y}^{i}(t)+g_{Z}(t)\hat{Z}^{i}(t)+g_{v_{i}}(t)v_{i}(t),\\
\hat{\bar{g}}^{i}(t)=&\bar{g}_{Y}(t)\hat{Y}^{i}(t)+\bar{g}_{Z}(t)\hat{Z}^{i}(t)+\bar{g}_{v_{i}}(t)v_{i}(t),\\
\hat{f}^{i}(t)=&f_{y}(t)\hat{y}^{i}(t)+f_{z}(t)\hat{z}^{i}(t)+f_{Y}(t)\hat{Y}^{i}(t)+f_{Z}(t)\hat{Z}^{i}(t)+f_{v_{i}}(t)v_{i}(t),\\
\hat{\bar{f}}^{i}(t)=&\bar{f}_{y}(t)\hat{y}^{i}(t)+\bar{f}_{z}(t)\hat{z}^{i}(t)+\bar{f}_{Y}(t)\hat{Y}^{i}(t)+\bar{f}_{Z}(t)\hat{Z}^{i}(t)+\bar{f}_{v_{i}}(t)v_{i}(t),\\
\hat{l}^{i}(t)=&l_{iy}(t)\hat{y}^{i}(t)+l_{iz}(t)\hat{z}^{i}(t)+l_{iY}(t)\hat{Y}^{i}(t)+l_{iZ}(t)\hat{Z}^{i}(t)+l_{iv_{i}}(t)v_{i}(t).\\
\end{align*}
Next, we define the generalized {\it
Hamiltonian function}
$H_{i}:[0,T]\times{\mathrm{R}}^n\times{\mathrm{R}}^{n\times
l}\times{\mathrm{R}}^m\times{\mathrm{R}}^{m\times
d}\times{\mathrm{U_{1}}}\times{\mathrm{U_{2}}}\times{\mathrm{R}}^n\times{\mathrm{R}}^{n\times
l}\times{\mathrm{R}}^m\times{\mathrm{R}}^{m\times d}$ as follows:
\begin{align}\label{DG15}
 H_{i}(t, y, z,& Y, Z, v_{1},v_{2},p_{i}, \bar{p}_{i}, q_{i}, \bar{q}_{i})
\triangleq\mbox{\ }\langle
      q_{i},f(y, z, Y, Z, v_{1},v_{2})\rangle+\langle \bar{q_{i}},\bar{f}(y, z, Y, Z,
v_{1},v_{2})\rangle\nonumber\\
&-\langle p_{i},g(Y, Z, v_{1}, v_{2})\rangle \mbox{\ }-\langle
\bar{p}_{i}, \bar{g}(Y, Z, v_{1}, v_{2})\rangle \mbox{\
}+l_{i}(y, z, Y, Z, v_{1}, v_{2}).
\end{align}
Let $(u_{1}, u_{2})\in
\mathcal{U}_{1}\times\mathcal{U}_{2}$ with the solution
$\big(y(\cdot),
z(\cdot), Y(\cdot), Z(\cdot)\big)$ of equation \eqref{DG14}.  We shall use the abbreviated notation $H_{i}(t)$
defined by
$$H_{i}(t)\equiv
H_{i}\big(t,y(t), z(t), Y(t), Z(t), u_{1}(t), u_{2}(t), p_{i}(t),
\bar{p}_{i}(t), q_{i}(t), \bar{q}_{i}(t)\big).$$ The adjoint
equations are described by the following generalized  stochastic
Hamiltonian systems:
\begin{equation}\label{DG22}
\left\{
\begin{aligned}
 dp_{i}(t)=&-H_{iY}^{ *}(t)dt-H_{iZ}^{ *}(t)dW(t)-\bar{p}_{i}(t)\hat{d}B(t),\\
-dq_{i}(t)=&\mbox{\ }H_{iy}^{*}(t)dt+H_{iz}^{ *}(t)\hat{d}B(t)-\bar{q_{i}}(t)dW(t),\\
  p_{i}(0)=& -\varphi_{iY}^{*}(Y(0))-\phi^{*}_{Y}\big(Y(0)\big)q_{i}(0),\\
  q_{i}(T)=&\gamma_{iy}^{*}\big(y(T)\big).
\end{aligned}
\right.
\end{equation}
Then we have the following maximum principle for nonzero-sum
differential games.
\begin{theorem}[Maximum principle for nonzero-sum games]\label{Theorem2.1}
Let (H1) hold and $\Big(u_{1}(\cdot),$ $u_{2}(\cdot)\Big)$ be an
equilibrium point of {\it Problem (NZSG)} with the corresponding
solutions $\Big(x(\cdot),y(\cdot),$ $z(\cdot)\Big)$ and $
\Big(p_{i}(\cdot), q_{i}(\cdot), k_{i}(\cdot)\Big)$  of \eqref{DG14}
and \eqref{DG22}.
 Then it follows that
\begin{equation}\label{DG5}
\begin{aligned}
\Big\langle E[ H_{1v_{1}}^{*}(t)|\mathcal{E}_{t}], v_{1}(t)-u_{1}(t)\Big\rangle\leq 0
\end{aligned}
\end{equation}
and
\begin{equation}\label{DG4}
\begin{aligned}
\Big\langle E[ H_{2v_{2}}^{*}(t)|\mathcal{E}_{t}], v_{2}(t)-u_{2}(t)\Big\rangle\leq 0
\end{aligned}
\end{equation}
are true for any $(v_{1}(\cdot),
v_{2}(\cdot))\in\mathcal{U}_{1}\times\mathcal{U}_{2}, a.e.\ a.s. $
\end{theorem}
{\it Proof:} Since $(u_{1}(\cdot), u_{2}(\cdot))$ is an
equilibrium point, we have  $$\frac{\partial h_{1}}{\partial \epsilon}(0, 0)=\lim\limits_{\epsilon\rightarrow 0}\frac{J_{1}(u_{1}+\epsilon v_{1}, u_{2})-J_{1}(u_{1}, u_{2})}{\epsilon}\leq0.$$
Then
\begin{align}
    0\geq&\frac{\partial}{\partial \epsilon}h_{1}(\epsilon, 0)|_{\epsilon=0}\nonumber\\
     =&
     E\int_{0}^{T}\Big(l_{1y}(t)\hat{y}^{1}(t)+l_{1z}(t)\hat{z}^{1}(t)+l_{1Y}(t)\hat{Y}^{1}(t)+l_{1Z}(t)\hat{Z}^{1}(t)+l_{1v_{1}}(t)v_{1}(t)\Big)dt\nonumber\\
     &+E\Big(\varphi_{1Y}\big(Y(0)\big)\hat{Y}^{1}(0)+\gamma_{1y}\big(y(T)\big)\hat{y}^{1}(T)\Big).\label{DG9}
\end{align}
Applying It\^{o}'s formula to $\langle p_{1}(t),
\hat{Y}^{1}(t)\rangle$ and $\langle q_{1}(t),
\hat{y}^{1}(t)\rangle$, and integrating from 0 to $T$, we have
\begin{align}
    &E\Big(\varphi_{1Y}\big(Y(0)\big)\hat{Y}^{1}(0)\Big)\nonumber\\
=& -E\langle p_{1}(0)+\phi_{Y}^{*}(Y(0))q_{1}(0),
    \hat{Y}^{1}(0)\rangle=-\langle\phi_{Y}^{*}(Y(0))q_{1}(0), \hat{Y}^{1}(0)\rangle\nonumber\\
&-E\int_{0}^{T}\Big(p_{1}^{*}(t)g_{v_{1}}(t)v_{1}(t)+q_{1}^{*}(t)f_{Y}(t)\hat{Y}^{1}(t)+\bar{q}_{1}^{*}(t)\bar{f}_{Y}(t)\hat{Y}^{1}(t)-l_{1Y}(t)\hat{Y}^{1}(t)\nonumber\\
&\hspace{1.3cm}+\bar{p}_{1}^{*}(t)\bar{g}_{v_{1}}(t)v_{1}(t)+q_{1}^{*}(t)f_{Z}(t)\hat{Z}^{1}(t)+\bar{q}_{1}^{*}(t)\bar{f}_{Z}(t)\hat{Z}^{1}(t)-l_{1Z}(t)\hat{Z}^{1}(t)\Big)dt,\label{DG8}
\end{align}
and
\begin{align}
    &E\Big(\gamma_{1y}\big(y(T)\big)\hat{y}^{1}(T)\Big)=\langle\phi_{Y}^{*}(Y(0))q_{1}(0),
    \hat{Y}^{1}(0)\rangle\nonumber\\
&+E\int_{0}^{T}\Big(q_{1}^{*}(t)f_{Y}(t)\hat{Y}^{1}(t)+q_{1}^{*}(t)f_{Z}(t)\hat{Z}^{1}(t)+q_{1}^{*}(t)f_{v_{1}}(t)v_{1}(t)-l_{1y}(t)\hat{y}^{1}(t)\nonumber\\
&\hspace{1.3cm}-l_{1z}(t)\hat{z}^{1}(t)+\bar{q}_{1}^{*}(t)\bar{f}_{Y}(t)\hat{Y}^{1}(t)+\bar{q}_{1}^{*}(t)\bar{f}_{Z}(t)\hat{Z}^{1}(t)+\bar{q}_{1}^{*}(t)\bar{f}_{v_{1}}(t)v_{1}(t)\Big)dt.\label{DG7}
\end{align}
Substituting \eqref{DG8} and
\eqref{DG7} into \eqref{DG9}, for all $v_{1}\in U_{1}$ such that
$u_{1}(\cdot)+v_{1}(\cdot)\in\mathcal{U}_{1},$ we get
\begin{align}
    0\geq&\frac{\partial}{\partial \epsilon}h_{1}(\epsilon, 0)|_{\epsilon=0}\nonumber\\
     =& E\int_{0}^{T}\Big(q_{1}^{*}(t)f_{v_{1}}(t)+\bar{q}_{1}^{*}(t)\bar{f}_{v_{1}}(t)+p_{1}^{*}(t)g_{v_{1}}(t)+\bar{p}_{1}^{*}(t)\bar{g}_{v_{1}}(t)+l_{1v_{1}}(t)\Big)v_{1}(t)dt\nonumber\\
     =& E\int_{0}^{T}\Big\langle H^{*}_{1v_{1}}(t), v_{1}(t)\Big\rangle dt=E\int_{0}^{T}E\left[\Big\langle H^{*}_{1v_{1}}(t), v_{1}(t)\Big\rangle\Big|\mathcal{E}_{t}\right] dt,\label{DG6}
\end{align}
 which implies that \eqref{DG5} is true.
The result \eqref{DG4} can be proved by the same method as shown in proving \eqref{DG5}.\hfill$\Box$

If the control process $\big(v_{1}(\cdot), v_{2}(\cdot))$ is admissible adapted to the filtration $\mathcal{F}_{t}$, we have the
following corollary.
\begin{corollary}[Maximum principle for full information nonzero-sum games]
Suppose that $\mathcal{E}_{t}=\mathcal{F}_{t}$ for all $t$.
Let (H1) hold and $\Big(u_{1}(\cdot),$ $u_{2}(\cdot)\Big)$ be an
equilibrium point of nonzero-sum differential games with the corresponding
solutions $\Big(x(\cdot),y(\cdot),$ $z(\cdot)\Big)$ and $
\Big(p_{i}(\cdot), q_{i}(\cdot), k_{i}(\cdot)\Big)$  of \eqref{DG14}
and \eqref{DG22}.
 Then it follows that
\begin{equation*}
\begin{aligned}
\Big\langle  H_{1v_{1}}^{*}(t), v_{1}(t)-u_{1}(t)\Big\rangle\leq 0
\end{aligned}
\end{equation*}
and
\begin{equation*}
\begin{aligned}
\Big\langle  H_{2v_{2}}^{*}(t), v_{2}(t)-u_{2}(t)\Big\rangle\leq 0
\end{aligned}
\end{equation*}
are true for any $(v_{1}(\cdot),
v_{2}(\cdot))\in\mathcal{U}_{1}\times\mathcal{U}_{2}, a.e.\ a.s. $
\end{corollary}

In what follows, we proceed to establish a verification theorem, also called a sufficient condition,
for an equilibrium point. For this, we introduce an additional  condition as follows: \\
(H2) $\phi(Y)=MY$ where $M$ is a non-zero constant matrix with order
$n\times m$. $\varphi_{i} $ and $\gamma_{i}$ are concave in $Y$ and
$y$ $(i=1, 2)$, respectively.

\begin{theorem}[Verification theorem for nonzero-sum games]\label{Thm3.2}Let (H1) and
(H2) hold. Let  $(u_{1}(\cdot),
u_{2}(\cdot))\in\mathcal{U}_{1}\times\mathcal{U}_{2}$ be with the
corresponding solutions $(y, z, Y, Z)$ and $ (p_{i}, \bar{p}_{i},
q_{i}, \bar{q_{i}})$ of equations \eqref{DG14} and \eqref{DG22}.
Suppose
\begin{equation}\label{Eq73}
    \begin{aligned}
&\hat{H}_{1}\big(t, a, b, c, d)= \ \sup_{v_{1}\in U_{1}}H_{1}\big(t, a, b, c, d, v_{1}, u_{2}(t), p_{1}(t),  \bar{p}_{1}(t), q_{1}(t), \bar{q}_{1}(t)\big),\\
&\hat{H}_{2}\big(t,  a, b, c, d)= \ \sup_{v_{2}\in U_{2}}H_{2}\big(t, a, b, c, d, u_{1}(t), v_{2}, p_{2}(t),  \bar{p}_{2}(t), q_{2}(t), \bar{q}_{2}(t)\big)\\
 &\hbox{exist for all}\ (t,  a, b, c, d)\in [0,
T]\times\mathrm{R}^{n}\times\mathrm{R}^{n\times
l}\times\mathrm{R}^{m}\times\mathrm{R}^{m\times
d}, \hbox{and are}\\
&\hbox{concave in}\ ( a, b, c, d)\ \hbox{for all t}\in [0, T]\
\hbox{(the Arrow condition)}.
    \end{aligned}
\end{equation}
Moreover
\begin{align}
&\mathbb{E}\Big[H_{1}\big(t, \ y(t), z(t),  Y(t), Z(t), u_{1}(t),
u_{2}(t), p_{1}(t), \bar{ p}_{1}(t), q_{1}(t),
\bar{q}_{1}(t)\big)\big|\mathcal{E}_{t}\Big]\nonumber\\
&=\sup_{v_{1}\in\ U_{1}} \mathbb{E}\Big[H_{1}\big(t,  y(t), z(t),
Y(t), Z(t), v_{1}, u_{2}(t), p_{1}(t), q_{1}(t),
\bar{q}_{1}(t)\big)\big|\mathcal{E}_{t}\Big],\label{Eq60}\\
&\mathbb{E}\Big[H_{2}\big(t, \ y(t), z(t),  Y(t), Z(t), u_{1}(t),
u_{2}(t),
p_{2}(t), \bar{p}_{2}(t), q_{2}(t), \bar{q}_{2}(t)\big)\big|\mathcal{E}_{t}\Big]\nonumber\\
&=\sup_{v_{2}\in\ U_{2}} \mathbb{E}\Big[H_{2}\big(t,  y(t), z(t),
Y(t), Z(t), u_{1}(t), v_{2}, p_{2}(t),  \bar{p}_{2}(t), q_{2}(t),
\bar{q}_{2}(t)\big)\big|\mathcal{E}_{t}\Big].\label{Eq60'}
\end{align}
Then  $(u_{1}(\cdot), u_{2}(\cdot))$ is an equilibrium point of {\it Problem (NZSG)}.
\end{theorem}
{\it Proof}\ : Let
 $(v_{1}(\cdot), u_{2}(\cdot))$ and $(u_{1}(\cdot), v_{2}(\cdot))$
$\in\mathcal{U}_{1}\times\mathcal{U}_{2}$ with the corresponding
solutions $(y^{v_{1}}, z^{v_{1}},$ $ Y^{v_{1}}, Z^{v_{1}})$ and $(
y^{v_{2}}, z^{v_{2}}, Y^{v_{2}}, Z^{v_{2}})$ to equation
\eqref{DG14}. We define the following terms
\begin{equation*}
    \begin{aligned}
    &H_{1}(t)=H_{1}(t,  y(t), z(t),  Y(t), Z(t), u_{1}(t), u_{2}(t), p_{1}(t), \bar{ p}_{1}(t), q_{1}(t),
\bar{q}_{1}(t)),\\
    &H_{1}^{v_{1}}(t)=H_{1}(t, y^{v_{1}}(t), z^{v_{1}}(t),  Y^{v_{1}}(t), Z^{v_{1}}(t), v_{1}(t), u_{2}(t), p_{1}(t), \bar{ p}_{1}(t), q_{1}(t),
\bar{q}_{1}(t)),\\
    &H_{1}^{v_{2}}(t)=H_{1}(t, y^{v_{2}}(t), z^{v_{2}}(t),  Y^{v_{2}}(t), Z^{v_{2}}(t), u_{1}(t), v_{2}(t), p_{1}(t), \bar{ p}_{1}(t), q_{1}(t),
\bar{q}_{1}(t)),\\
    &f^{v_{1}}(t)=f(t, y^{v_{1}}(t), z^{v_{1}}(t),  Y^{v_{1}}(t), Z^{v_{1}}(t),  v_{1}(t), u_{2}(t)),\\
    & f^{v_{2}}(t)=f(t,y^{v_{2}}, z^{v_{2}},  Y^{v_{2}}, Z^{v_{2}}, u_{1}(t), v_{2}(t)),
    \end{aligned}
\end{equation*}
and similar notations are made for $\bar{f}^{v_{1}},
\bar{f}^{v_{2}},\cdots$.

By virtue of the concavity property of $\varphi_{1} $ and
$\gamma_{1}$, we have for $\forall\ v_{1}(\cdot)\in \mathcal{U}_{1}$
\begin{equation}\label{Sq1}
    \begin{aligned}
    J_{1}(v_{1}(\cdot), u_{2}(\cdot))-J_{1}(u_{1}(\cdot), u_{2}(\cdot))\leq
    I_{1}+I_{2}+I_{3}
    \end{aligned}
\end{equation}
with
\begin{eqnarray*}
    \begin{aligned}
I_{1}&=\mathbb{E}\left[ \gamma_{1y}(y(T))(y^{v_{1}}(T)-y(T))\right],\\
I_{2}&=\mathbb{E}\left[ \varphi_{1Y}(Y(0))(Y^{v_{1}}(0)-Y(0))\right],\\
I_{3}&=\mathbb{E}\int_{0}^{T}\Big(l^{v_{1}}_{1}(t)-l_{1}(t)\Big)dt.
   \end{aligned}
\end{eqnarray*}
Applying It\^{o}'s formula to $ \langle q_{1}(t),
y^{v_{1}}(t)-y(t)\rangle$ and $ \langle p_{1}(t),
Y^{v_{1}}(t)-Y(t)\rangle,$
\begin{align}
    I_{1} =&\, \mathbb{E}[\langle q_{1}(0), M(Y^{v_{1}}(0)-Y(0))\rangle]\nonumber\\
     &+\mathbb{E} \int_{0}^{T}\Big(\langle q_{1}(t), f^{v_{1}}(t)-f(t)\rangle-\langle H_{1y}^{*}(t), y^{v_{1}}(t)-y(t)\rangle\nonumber\\
     &\hspace{2cm}+\langle \bar{q}_{1}(t), \bar{f}^{v_{1}}(t)-\bar{f}(t)\rangle-\langle H_{1z}^{*}(t), z^{v_{1}}(t)-z(t)\rangle\Big)dt,\label{Eq23'}\\
  I_{2}= &\,-\mathbb{E}[\langle q_{1}(0), M(Y^{v_{1}}(0)-Y(0))\rangle]\nonumber\\
         &-\mathbb{E}\int_{0}^{T}\Big(\langle p_{1}(t), g^{v_{1}}(t)-g(t)\rangle+\langle  H_{1Y}^{*}(t), Y^{v_{1}}(t)-Y(t)\rangle\nonumber\\
       &+\langle \bar{p}_{1}(t), \bar{g}^{v_{1}}(t)-\bar{g}(t)\rangle+\langle H_{1Z}^{*}(t), Z^{v_{1}}(t)-Z(t)\rangle
       \Big)dt,\\
  I_{3}=\ & \mathbb{E}\int_{0}^{T}\bigg(H_{1}^{v_{1}}(t)-H_{1}(t)-\langle q_{1}(t), f^{v_{1}}(t)-f(t)\rangle-\langle \bar{q}_{1}(t), \bar{f}^{v_{1}}(t)-\bar{f}(t)\rangle\nonumber\\
       &\hspace{1.7cm}+\langle \bar{p}_{1}(t), \bar{g}^{v_{1}}(t)-\bar{g}(t)\rangle
        +\langle p_{1}(t), g^{v_{1}}(t)-g(t)\rangle
        \bigg)dt.\label{Eq23}
\end{align}
Substituting \eqref{Eq23'}---\eqref{Eq23} into \eqref{Sq1}, it
follows immediately that
\begin{align}
   & J_{1}(v_{1}(\cdot), u_{2}(\cdot))-J_{1}(u_{1}(\cdot), u_{2}(\cdot))\nonumber\\
 \leq&\ \mathbb{E}\int_{0}^{T}\bigg(H_{1}^{v_{1}}(t)-H_{1}(t)-\langle
H_{1Y}^{*}(t), Y^{v_{1}}(t)-Y(t)\rangle-\langle
H_{1Z}^{*}(t), Z^{v_{1}}(t)-Z(t)\rangle\nonumber\\
&\hspace{1.5cm} -\langle H_{1y}^{*}(t),
y^{v_{1}}(t)-y(t)\rangle-\langle H_{1z}^{*}(t),
z^{v_{1}}(t)-z(t)\rangle\bigg) dt.\label{Eq62}
   \end{align}
Since $v_{1}\longrightarrow\mathbb{E}\Big[H_{1}\big(t,  y(t), z(t),
Y(t), Z(t), v_{1}, u_{2}(t), p_{1}(t), q_{1}(t),
\bar{q}_{1}(t)\big)\big|\mathcal{E}_{t}\Big]$ is maximum for
$v_{1}=u_{1}$ and since $v_{1}(t), u_{1}(t)$ are
$\mathcal{E}_{t}$-measurable, we get
\begin{align*}
    &\mathbb{E}\Big[\frac{\partial}{\partial v_{1}}H_{1}\big(t, \ y(t), z(t),  Y(t), Z(t), u_{1}(t),
u_{2}(t), p_{1}(t), \bar{ p}_{1}(t), q_{1}(t),
\bar{q}_{1}(t)\big)\big(v_{1}(t)-u_{1}(t)\big)\big|\mathcal{E}_{t}\Big]\\
=&\mathbb{E}\Big[\frac{\partial}{\partial v_{1}}H_{1}\big(t, \ y(t),
z(t),  Y(t), Z(t), v_{1}(t), u_{2}(t), p_{1}(t), \bar{ p}_{1}(t),
q_{1}(t),
\bar{q}_{1}(t)\big)\big|\mathcal{E}_{t}\Big]_{v_{1}=u_{1}}\big(v_{1}(t)-u_{1}(t)\big)\\
\leq& 0.
\end{align*}
By the equality \eqref{Eq60} and the concavity of $\hat{H}_{1}$, we
conclude that
\begin{equation}\label{Eq61}
\begin{aligned}
J_{1}(v_{1}(\cdot), u_{2}(\cdot))-J_{1}(u_{1}(\cdot),
u_{2}(\cdot))\leq 0,
\end{aligned}
\end{equation}
for all $v_{1}(\cdot)\in \mathcal{U}_{1}$.
Repeating the similar proceeding as shown in deriving \eqref{Eq61},
we can prove that
\begin{equation}\label{Eq35'}
    J_{2}(u_{1}(\cdot), v_{2}(\cdot))-J_{2}(u_{1}(\cdot), u_{2}(\cdot))\leq 0.
\end{equation}
Based on the arguments above,  $(u_{1}(\cdot), u_{2}(\cdot))$ is an
equilibrium point of {\it Problem (NZSG)}.\hfill $\Box$

\begin{corollary}[Verification theorem for full information nonzero-sum games]Suppose that $\mathcal{E}_{t}=\mathcal{F}_{t}$ for all $t$ and that (H1),
(H2) and \eqref{Eq73} hold. Let  $(u_{1}(\cdot),
u_{2}(\cdot))\in\mathcal{U}_{1}\times\mathcal{U}_{2}$ be with the
corresponding solutions $(y, z, Y, Z)$ and $ (p_{i}, \bar{p}_{i},
q_{i}, \bar{q_{i}})$ of equations \eqref{DG14} and \eqref{DG22}.
Moreover
\begin{align*}
&H_{1}\big(t, \ y(t), z(t),  Y(t), Z(t), u_{1}(t),
u_{2}(t), p_{1}(t), \bar{ p}_{1}(t), q_{1}(t),
\bar{q}_{1}(t)\big)\nonumber\\
&=\sup_{v_{1}\in\ U_{1}} H_{1}\big(t,  y(t), z(t),
Y(t), Z(t), v_{1}, u_{2}(t), p_{1}(t), q_{1}(t),
\bar{q}_{1}(t)\big)\\
&H_{2}\big(t, \ y(t), z(t),  Y(t), Z(t), u_{1}(t),
u_{2}(t),
p_{2}(t), \bar{p}_{2}(t), q_{2}(t), \bar{q}_{2}(t)\big)\nonumber\\
&=\sup_{v_{2}\in\ U_{2}} H_{2}\big(t,  y(t), z(t),
Y(t), Z(t), u_{1}(t), v_{2}, p_{2}(t),  \bar{p}_{2}(t), q_{2}(t),
\bar{q}_{2}(t)\big).
\end{align*}
Then  $(u_{1}(\cdot), u_{2}(\cdot))$ is an equilibrium point of nonzero-sum differential games.
\end{corollary}

\subsection{Zero-sum case }

In this section, we consider zero-sum differential games of FBDSDEs.
In fact,  zero-sum games can be consider  a special case of nonzero-sum games.
By the maximum principle of nonzero-sum games in Section 3.1,
we can deduce the necessary conditions for a saddle point of zero-sum games. We shall detail this as follows.

Let
$$-J_{1}=J_{2}=J.$$
If $(u_{1}(\cdot), u_{2}(\cdot))$ is an equilibrium point of
{\it Problem (NZSG)}, we have
\begin{equation*}
    \left\{
\begin{aligned}
&J_{1}(u_{1}(\cdot), u_{2}(\cdot))\geq J_{1}(v_{1}(\cdot),
u_{2}(\cdot)),\\
&J_{2}(u_{1}(\cdot), u_{2}(\cdot))\geq J_{2}(u_{1}(\cdot),
v_{2}(\cdot)),
\end{aligned}
    \right.
\end{equation*}
which implies that
\begin{equation*}
\begin{aligned}
 J(u_{1}(\cdot), v_{2}(\cdot))\leq J(u_{1}(\cdot), u_{2}(\cdot))\leq J(v_{1}(\cdot),
 u_{2}(\cdot)).
\end{aligned}
\end{equation*}
Therefore, $(u_{1}(\cdot), u_{2}(\cdot))$ is also a saddle
point of {\it Problem (ZSG)}.

We define a new {\it
Hamiltonian function}
$H:[0,T]\times{\mathrm{R}}^n\times{\mathrm{R}}^{n\times
l}\times{\mathrm{R}}^m\times{\mathrm{R}}^{m\times
d}\times{\mathrm{U_{1}}}\times{\mathrm{U_{2}}}\times{\mathrm{R}}^n\times{\mathrm{R}}^{n\times
l}\times{\mathrm{R}}^m\times{\mathrm{R}}^{m\times d}$ as follows:
\begin{align}\label{DG3}
 H(t, y, z,& Y, Z, v_{1},v_{2}, p, \bar{p}, q, \bar{q})
\triangleq\mbox{\ }\langle
      q, f(y, z, Y, Z, v_{1},v_{2})\rangle+\langle \bar{q},\bar{f}(y, z, Y, Z,
v_{1},v_{2})\rangle\nonumber\\
&-\langle p, g(Y, Z, v_{1}, v_{2})\rangle \mbox{\ }-\langle
\bar{p}, \bar{g}(Y, Z,  v_{1}, v_{2})\rangle \mbox{\
}+l(y, z, Y, Z, v_{1}, v_{2}).
\end{align}
Let $(u_{1}, u_{2})\in
\mathcal{U}_{1}\times\mathcal{U}_{2}$ with the solution
$\big(y(\cdot),
z(\cdot), Y(\cdot), Z(\cdot)\big)$ of equation \eqref{DG14}.  We shall use the abbreviated notation $H(t)$
defined by
$$H(t)\equiv
H\big(t,y(t), z(t), Y(t), Z(t), u_{1}(t), u_{2}(t), p(t),
\bar{p}(t), q(t), \bar{q}(t)\big).$$ The adjoint
equations are described by the following generalized  stochastic
Hamiltonian systems:
\begin{equation}\label{DG2}
\left\{
\begin{aligned}
 dp(t)=&-H_{Y}^{ *}(t)dt-H_{Z}^{ *}(t)dW(t)-\bar{p}(t)\hat{d}B(t),\\
-dq(t)=&\mbox{\ }H_{y}^{*}(t)dt+H_{z}^{ *}(t)\hat{d}B(t)-\bar{q}(t)dW(t),\\
  p(0)=& -\varphi_{Y}^{*}(Y(0))-\phi^{*}_{Y}\big(Y(0)\big)q(0),\\
  q(T)=&\gamma_{y}^{*}\big(y(T)\big).
\end{aligned}
\right.
\end{equation}

Based on the above arguments, we can directly  derive the following
maximum principle for zero-sum games.
\begin{theorem}[Maximum principle for zero-sum games]\label{Theorem2.3}
Let (H1) hold and $(u_{1}(\cdot), u_{2}(\cdot))$ be a saddle point of {\it Problem (ZSG)} with the solutions   $\Big(y(\cdot), z(\cdot), Y(\cdot), Z(\cdot)\Big)$ and $
\Big(p(\cdot), \bar{p}(\cdot), q(\cdot),$ $\bar{q}(\cdot)\Big)$ to \eqref{DG14}
and \eqref{DG2}, respectively.
 Then it follows that
\begin{equation}
\begin{aligned}
\Big\langle E[ H_{v_{1}}^{*}(t)|\mathcal{E}_{t}], v_{1}(t)-u_{1}(t)\Big\rangle\geq 0
\end{aligned}
\end{equation}
and
\begin{equation}
\begin{aligned}
\Big\langle E[ H_{v_{2}}^{*}(t)|\mathcal{E}_{t}], v_{2}(t)-u_{2}(t)\Big\rangle\leq 0
\end{aligned}
\end{equation}
are true for any $(v_{1}(\cdot),
v_{2}(\cdot))\in\mathcal{U}_{1}\times\mathcal{U}_{2}, a.e.\ a.s. $
\end{theorem}
\begin{remark}
If $\big(u_{1}(\cdot), u_{2}(\cdot)\big)$ is an equilibrium point (resp. a saddle point) of nonzero-sum (resp. zero-sum) differential games and $\big(u_{1}(t), u_{2}(t)\big)$ is an interior point of $U_{1}\times U_{2}$ a.s. for all $t\in [0, T]$, then the inequalities in Theorem 3.1 (resp. Theorem 3.3) are equivalent to the following equations
\begin{align*}
  E[ H_{iv_{i}}^{*}(t)|\mathcal{E}_{t}]=0, i=1, 2\  \big(resp.\ E[ H_{v_{j}}^{*}(t)|\mathcal{E}_{t}]=0, j=1,2\big).
\end{align*}

\end{remark}

We note that Theorem \ref{Theorem2.3} gives a globally necessary condition for a saddle point of zero-sum games. In the following, we begin to present a corresponding locally necessary condition for {\it Problem (ZSG)}.

We firstly give the following assumptions.\\
(H3) For all $t, \tau$ such that $0\leq t<t+\tau\leq T,$ all bounded
$\mathcal{E}_{t}$-measurable $\alpha_{1}, \alpha_{2},$ and for $s\in
[0, T]$, the control $\beta_{1}(s)\doteq (0, \cdots, \beta_{1j}(s),
\cdots, 0) $ and $\beta_{2}(s)\doteq (0, \cdots, \beta_{2j}(s),
\cdots, 0) $, $j=1, \cdots, n$, with
$$\beta_{1j}(s)\doteq \alpha_{1j}\chi_{[t,
t+\tau]}(s)\;\;\;\hbox{and}\;\;\;
   \beta_{2j}(s)\doteq \alpha_{2j}\chi_{[t, t+\tau]}(s)$$
belong to $\mathcal{U}_{1}$ and $\mathcal{U}_{2}$. For given $u_{1},
\beta_{1}\in \mathcal{U}_{1}$ and $u_{2}, \beta_{2}\in
\mathcal{U}_{2}$, there exists $\delta>0$ such that
$$u_{1}+\epsilon\beta_{1}\in\mathcal{U}_{1}\quad \hbox{and} \quad u_{2}+\rho\beta_{2}\in\mathcal{U}_{1},$$
where $\beta_{1}$ and $\beta_{2}$ are bounded, and $\epsilon,\rho
\in (-\delta, \delta).$

\begin{theorem}[Local maximum principle for zero-sum games]
Let (H1) and (H3) hold. Let $\big(u_{1}(\cdot),
u_{2}(\cdot)\big)\in\mathcal{U}_{1}\times\mathcal{U}_{2}$ with the
solutions $\big(y(\cdot), z(\cdot), Y(\cdot), Z(\cdot)\big)$ and
$\big(p(\cdot), \bar{p}(\cdot), q(\cdot),
\bar{q}(\cdot)\big)$ to equations \eqref{DG14} and \eqref{DG2},
respectively. Further, $\big(u_{1}(\cdot), u_{2}(\cdot)\big)$ is a
directional critical point for $J\big(v_{1}(\cdot),
v_{2}(\cdot)\big)$, in the sense that, for all bounded $\beta_{1}\in
\mathcal{U}_{1}$ and $\beta_{2}\in \mathcal{U}_{2}$, there exists
$\delta>0$ such that $u_{1}+\epsilon\beta_{1}\in\mathcal{U}_{1}\quad
\hbox{and} \quad u_{2}+\rho\beta_{2}\in\mathcal{U}_{2},$  and
$$h(\epsilon, \rho)\doteq J(u_{1}+\epsilon\beta_{1},
u_{2}+\rho\beta_{2})$$ has a critical point at (0, 0), for all
$\epsilon, \rho\in(-\delta, \delta)$, i.e.
$$\frac{\partial h}{\partial \epsilon}(0, 0)=\frac{\partial h}{\partial \rho}(0, 0)=0.$$
\end{theorem}
Then, for  all $t\in[0, T]$, we have
\begin{align}\label{DG1}
   E[H_{v_{1}}(t)|\mathcal{E}_{t}]=E[H_{v_{2}}(t)|\mathcal{E}_{t}]=0.
\end{align}
{\it Proof:} Since the remaining case can be dealt with by the similar
proceeding, we only prove
\begin{align*}
    E[H_{v_{1}}(t)|\mathcal{E}_{t}]=0.
\end{align*}
For $\frac{\partial h}{\partial \epsilon}(0, 0)=0,$ we have
\begin{align}
    0=&\frac{\partial}{\partial \epsilon}h(\epsilon, 0)|_{\epsilon=0}\nonumber\\
     =&
     E\int_{0}^{T}\Big(l_{y}(t)\hat{y}^{1}(t)+l_{z}(t)\hat{z}^{1}(t)+l_{Y}(t)\hat{Y}^{1}(t)+l_{Z}(t)\hat{Z}^{1}(t)+l_{v_{1}}(t)\beta_{1}\Big)dt\nonumber\\
     &+E\Big(\varphi_{Y}\big(Y(0)\big)\hat{Y}^{1}(0)+\gamma_{y}\big(y(T)\big)\hat{y}^{1}(T)\Big).\label{DG23}
\end{align}
As shown in Theorem \ref{Theorem2.1}, applying It\^{o}'s formula to $\langle p(t),
\hat{Y}^{1}(t)\rangle$ and $\langle q(t),
\hat{y}^{1}(t)\rangle$, integrating from 0 to $T$, and substituting them into \eqref{DG23},  we have
\begin{align}
    0=&\frac{\partial}{\partial \epsilon}h(\epsilon, 0)|_{\epsilon=0}\nonumber\\
     =& E\int_{0}^{T}\Big(q^{*}(t)f_{v_{1}}(t)+\bar{q}^{*}(t)\bar{f}_{v_{1}}(t)+p^{*}(t)g_{v_{1}}(t)+\bar{p}^{*}(t)\bar{g}_{v_{1}}(t)+l_{v_{1}}(t)\Big)\beta_{1}(t)dt\nonumber\\
     =& E\int_{0}^{T}H^{*}_{v_{1}}(t)\beta_{1}(t)dt.\label{DG24}
\end{align}
In terms of assumption (H3) and equality \eqref{DG24}, we further derive that
\begin{align*}
   E\int_{t}^{t+\tau}H_{v_{1j}}(s)\beta_{1j}(s)ds=0.
\end{align*}
Differentiating with respect to $\tau$ at  $\tau=0$, it yields that
\begin{align}\label{DG25}
   E[H_{v_{1j}}(t)\beta_{1j}(t)]=0.
\end{align}
Since equality \eqref{DG25} holds for all bounded
$\mathcal{E}_{t}$-measurable $\beta_{1j}$, we conclude that
\begin{align*}
   E[H_{v_{1}}(t)|\mathcal{E}_{t}]=0.
\end{align*}
Repeating the similar proceeding by differentiating the function
$h(0, \rho)$ with respect to $\rho$, we have
\begin{align*}
   E[H_{v_{2}}(t)|\mathcal{E}_{t}]=0.
\end{align*}
The proof is completed.\hfill $\Box$

In the sequel, we give a verification theorem for a saddle point of  zero-sum
games.

\begin{theorem}[Verification theorem for zero-sum games] Let $(H_{1})$ hold and  $\phi(Y)=MY$ where $M$ is a non-zero constant matrix with order
$n\times m$. Let $(u_{1}(\cdot),
u_{2}(\cdot))\in\mathcal{U}_{1}\times\mathcal{U}_{2}$ with the
solutions $(y, z, Y, Z)$ and $(p, \bar{p}, q, \bar{q})$ to equations
\eqref{DG14} and \eqref{DG2}, respectively. Suppose that the {\it Hamiltonian
function H} satisfies the following conditional mini-maximum
principle:
\begin{align}
    &\ E\Big[ H\big(t,  y(t), z(t), Y(t), Z(t), u_{1}(t), u_{2}(t), p(t), \bar{p}(t), q(t), \bar{q}(t)\big)\big|\mathcal{E}_{t}\Big]\nonumber\\
    =\ &\inf_{v_{1}(\cdot)\in\mathcal{U}_{1}}E\Big[ H(t, y(t), z(t), Y(t), Z(t), v_{1}(t), u_{2}(t), p(t), \bar{p}(t), q(t), \bar{q}(t)\big)\big|\mathcal{E}_{t}\Big]\nonumber\\
    =\ &\sup_{v_{2}(\cdot)\in\mathcal{U}_{2}}E\Big[H(t, y(t), z(t), Y(t), Z(t), u_{1}(t), v_{2}(t),  p(t), \bar{p}(t), q(t), \bar{q}(t)\big)\big|\mathcal{E}_{t}\Big].
\end{align}
(i)Assume that both $\varphi$ and $\gamma$ are concave, and
\begin{equation*}
    \begin{aligned}
\hat{H}_{2}\big(t, a, b, c)=& \ \sup_{v_{2}(\cdot)\in \mathcal{U}_{2}}H\big(t, a, b, c, u_{1}(t), v_{2}, p(t), q(t), k(t)\big),\\
    \end{aligned}
\end{equation*}
exists for all $(t,  a, b, c)\in [0,
T]\times\mathrm{R}^{n}\times\mathrm{R}^{m}\times\mathrm{R}^{m\times
d}$, and is concave in $( a, b, c)$. Then we have
\begin{equation*}
     J(u_{1}(\cdot), v_{2}(\cdot))\leq J(u_{1}(\cdot), u_{2}(\cdot)),
    \quad\quad \hbox{for all}\; v_{2}(\cdot)\in \mathcal{U}_{2},
\end{equation*}
and
\begin{equation*}
    J(u_{1}(\cdot), u_{2}(\cdot))=\sup_{v_{2}(\cdot)\in \mathcal{U}_{2}}J(u_{1}(\cdot),
    v_{2}(\cdot)).
   \end{equation*}
(ii)Assume that both $\varphi$ and $\gamma$ are convex, and
\begin{equation*}
    \begin{aligned}
\hat{H}_{1}\big(t, a, b, c)=& \ \inf_{v_{1}(\cdot)\in \mathcal{U}_{1}}H\big(t, a, b, c, v_{1}, u_{2}(t), p(t), q(t), k(t)\big),\\
    \end{aligned}
\end{equation*}
exists for all $(t,  a, b, c)\in [0,
T]\times\mathrm{R}^{n}\times\mathrm{R}^{m}\times\mathrm{R}^{m\times
d}$, and is convex in $( a, b, c)$. Then we have
\begin{equation*}
    J(u_{1}(\cdot), u_{2}(\cdot))\leq J(v_{1}(\cdot), u_{2}(\cdot)),
    \quad\quad \hbox{for all}\; v_{1}(\cdot)\in \mathcal{U}_{1},
\end{equation*}
and
\begin{equation*}
    J(u_{1}(\cdot), u_{2}(\cdot))=\inf_{v_{1}(\cdot)\in \mathcal{U}_{1}}J(v_{1}(\cdot),
    u_{2}(\cdot)).
   \end{equation*}
(iii) If both (i) and (ii) are true, then $(u_{1}(\cdot),
u_{2}(\cdot))$ is a saddle point which implies
\begin{align*}
    \sup_{v_{2}(\cdot)\in \mathcal{U}_{2}}\Big(\inf_{v_{1}(\cdot)\in \mathcal{U}_{1}}J\big(v_{1}(\cdot),
    v_{2}(\cdot)\big)\Big)&=J(u_{1}(\cdot), u_{2}(\cdot))\\
                           &=\inf_{v_{1}(\cdot)\in \mathcal{U}_{1}}\Big(\sup_{v_{2}(\cdot)\in \mathcal{U}_{2}}J\big(v_{1}(\cdot),
    v_{2}(\cdot)\big)\Big).
\end{align*}

\end{theorem}
{\it Proof :} (i)  Using the same argument as in the proof of
Theorem \ref{Thm3.2}, we can the following:
\begin{equation}\label{Eq38'}
     J(u_{1}(\cdot), v_{2}(\cdot))\leq J(u_{1}(\cdot), u_{2}(\cdot)),
    \quad\quad \hbox{for all}\; v_{2}(\cdot)\in \mathcal{U}_{2}.
\end{equation}
Furthermore
\begin{equation*}
    \sup_{v_{2}(\cdot)\in \mathcal{U}_{2}} J(u_{1}(\cdot), v_{2}(\cdot))\leq J(u_{1}(\cdot),
    u_{2}(\cdot)).
\end{equation*}
 Since $u_{2}(\cdot)\in \mathcal{U}_{2}$,  we
have
\begin{equation*}
    \sup_{v_{2}(\cdot)\in \mathcal{U}_{2}} J(u_{1}(\cdot), v_{2}(\cdot))=J(u_{1}(\cdot),
    u_{2}(\cdot)).
\end{equation*}
(ii) This statement can be proved in a similar way as shown
before.\\
(iii) If both (i) and (ii) are true, then
\begin{equation*}
     J(u_{1}(\cdot), v_{2}(\cdot))\leq J(u_{1}(\cdot), u_{2}(\cdot))
     \leq  J(v_{1}(\cdot), u_{2}(\cdot)),
\end{equation*}
for  all $(v_{1}(\cdot), v_{2}(\cdot))\in
\mathcal{U}_{1}\times\mathcal{U}_{2}$, i.e. $(u_{1}(\cdot),
u_{2}(\cdot))$ is a saddle point.

In the following, on the one hand, we have
\begin{equation*}
    J(u_{1}(\cdot), u_{2}(\cdot))\leq\inf_{v_{1}(\cdot)\in \mathcal{U}_{1}} J(v_{1}(\cdot), u_{2}(\cdot))
    \leq \inf_{v_{1}(\cdot)\in \mathcal{U}_{1}}\Big(\sup_{v_{2}(\cdot)\in \mathcal{U}_{2}} J\big(v_{1}(\cdot)
    v_{2}(\cdot)\big)\Big),
\end{equation*}
and
\begin{equation*}
    J(u_{1}(\cdot), u_{2}(\cdot))\geq\sup_{v_{2}(\cdot)\in \mathcal{U}_{2}} J(u_{1}(\cdot), v_{2}(\cdot))
    \geq \sup_{v_{2}(\cdot)\in \mathcal{U}_{2}}\Big(\inf_{v_{1}(\cdot)\in \mathcal{U}_{1}} J\big(v_{1}(\cdot),
    v_{2}(\cdot)\big)\Big),
\end{equation*}
which imply that
\begin{align}\label{Eq70}
    \sup_{v_{2}(\cdot)\in \mathcal{U}_{2}}\Big(\inf_{v_{1}(\cdot)\in \mathcal{U}_{1}} J\big(v_{1}(\cdot),
    v_{2}(\cdot)\big)\Big)&\leq J(u_{1}(\cdot), u_{2}(\cdot))\nonumber\\
                          &\leq \inf_{v_{1}(\cdot)\in \mathcal{U}_{1}}\Big(\sup_{v_{2}(\cdot)\in \mathcal{U}_{2}} J\big(v_{1}(\cdot),
    v_{2}(\cdot)\big)\Big).
\end{align}
On the other hand, we have
\begin{equation*}
    J(u_{1}(\cdot), u_{2}(\cdot))\leq\inf_{v_{1}(\cdot)\in \mathcal{U}_{1}} J(v_{1}(\cdot), u_{2}(\cdot))
    \leq \sup_{v_{2}(\cdot)\in \mathcal{U}_{2}}\Big(\inf_{v_{1}(\cdot)\in \mathcal{U}_{1}}
    J\big(v_{1}(\cdot),
    v_{2}(\cdot)\big)\Big)
\end{equation*}
and
\begin{equation*}
    J(u_{1}(\cdot), u_{2}(\cdot))\geq\sup_{v_{2}(\cdot)\in \mathcal{U}_{2}} J(u_{1}(\cdot), v_{2}(\cdot))
    \geq \inf_{v_{1}(\cdot)\in \mathcal{U}_{1}}\Big(\sup_{v_{2}(\cdot)\in \mathcal{U}_{2}} J\big(v_{1}(\cdot),
    v_{2}(\cdot)\big)\Big),
\end{equation*}
which imply that
\begin{align}\label{Eq71}
    \sup_{v_{2}(\cdot)\in \mathcal{U}_{2}}\Big(\inf_{v_{1}(\cdot)\in \mathcal{U}_{1}} J\big(v_{1}(\cdot),
    v_{2}(\cdot)\big)\Big)&\geq J(u_{1}(\cdot), u_{2}(\cdot))\nonumber\\
                          &\geq \inf_{v_{1}(\cdot)\in \mathcal{U}_{1}}\Big(\sup_{v_{2}(\cdot)\in \mathcal{U}_{2}} J\big(v_{1}(\cdot),
    v_{2}(\cdot)\big)\Big).
\end{align}
Combining \eqref{Eq70} and \eqref{Eq71},  we have
\begin{align*}
    \sup_{v_{2}(\cdot)\in \mathcal{U}_{2}}\Big(\inf_{v_{1}(\cdot)\in \mathcal{U}_{1}}J\big(v_{1}(\cdot),
    v_{2}(\cdot)\big)\Big)&=J(u_{1}(\cdot), u_{2}(\cdot))\\
                           &=\inf_{v_{1}(\cdot)\in \mathcal{U}_{1}}\Big(\sup_{v_{2}(\cdot)\in
                           \mathcal{U}_{2}}J\big(v_{1}(\cdot),
    v_{2}(\cdot)\big)\Big).
\end{align*}
\begin{remark}
Similar to the results in Section 3.1, we can also give the corresponding corollaries for maximum principle and verification theorem of a saddle point of full information zero-sum differential games. For simplicity, we omit them here. \end{remark}

\section{An example on a nonzero-sum game }
In this section,  an example of nonzero-sum differential
games  of FBSDEs is worked out to illustrate our theoretical result. Firstly, by
applying the maximum principle (see Theorem 3.1), we
find a candidate equilibrium point. Then we obtain the explicit expression of the candidate equilibrium point by virtue of certain filtering techniques of forward-backward stochastic differential equations. Finally, using the verification theorem of an equilibrium point (see Theorem 3.2), we confirm that
it is indeed an equilibrium point.

{\bf Example:}
Consider the system of FBDSDE
\begin{equation}\label{DG19}
\left\{\begin{aligned}
-dY^{v_{1}, v_{2}}(t)=&\big[a_{0}(t)+a_{1}(t)Y^{v_{1}, v_{2}}(t)+a_{2}(t)Z^{v_{1}, v_{2}}(t)+F_{i2}(t)v_{1}(t)+a_{4}(t)v_{2}(t)\big]dt\\
                      &+b_{0}(t)\hat{d}B(t)-Z^{v_{1}, v_{2}}(t)dW(t),\\
dy^{v_{1}, v_{2}}(t)=&\Big[c_{0}(t)+c_{1}(t)y^{v_{1}, v_{2}}(t)+c_{2}(t)Y^{v_{1}, v_{2}}(t)+c_{3}(t)Z^{v_{1}, v_{2}}(t)\Big]dt\\
                     &+d_{0}(t)dW(t)-z^{v_{1}, v_{2}}(t)\hat{d}B(t),\\
Y^{v_{1}, v_{2}}(T)=&\ \xi,\;\;\;y^{v_{1}, v_{2}}(0)=\ MY^{v_{1}, v_{2}}(0),
\end{aligned}\right.
\end{equation}
with the performance criterion, for $i=1, 2,$
\begin{align}\label{DG20}
  J_{i}\big(&v_{1}(\cdot), v_{2}(\cdot)\big)
=-\frac{1}{2}\mathbb{E}\bigg[\int_{0}^{T}\Big(\langle e_{i1}(t)y^{v_{1}, v_{2}}(t), y^{v_{1}, v_{2}}(t)\rangle+\langle e_{i2}(t)z^{v_{1}, v_{2}}(t), z^{v_{1}, v_{2}}(t)\rangle\nonumber\\
   &+\langle e_{i3}(t)Y^{v_{1}, v_{2}}(t), Y^{v_{1}, v_{2}}(t)\rangle+\langle e_{i4}(t)Z^{v_{1}, v_{2}}(t), Z^{v_{1}, v_{2}}(t)\rangle+\langle e_{i7}(t)v_{i}(t), v_{i}(t)\rangle\Big)dt\nonumber\\
   &+\langle e_{i5}(T)y^{v_{1}, v_{2}}(T), y^{v_{1}, v_{2}}(T)\rangle+\langle e_{i6}Y^{v_{1}, v_{2}}(0), Y^{v_{1}, v_{2}}(0)\rangle\bigg].
\end{align}
Here, we assume that all the coefficients in \eqref{DG19} and \eqref{DG20}  are bounded and deterministic functions of $t$,
$e_{i1}, \cdots, e_{i6}$  are symmetric nonnegative definite,
 and  $e_{i7}$ is symmetric uniformly positive definite.
The set of admissible controls is defined by
\begin{equation}\label{DG29}
\begin{aligned}
    \mathcal{U}_{i}=\{v_{i}(\cdot)\;|
\;  v_{i}(\cdot)\ \hbox{is an}&\ \mathrm{R}^{k_{i}}\hbox{-valued}\
\mathcal{E}_{t}\hbox{-adapted process}\\
& \hbox{and satisfies}\ \mathbb{E}\int_{0}^{T}v_{i}^{2}(t)dt<\infty
\}, i=1, 2,
 \end{aligned}
\end{equation}
Where $$\mathcal{E}_{t}=\mathcal{N}\vee \sigma\big\{W(r): 0\leq r\leq t\big\}.$$
For simplicity, we only deal with the case of 1-dimensional coefficients.
Our problem is to find $(u_{1}(\cdot), u_{2}(\cdot))\in \mathcal{U}_{1}\times\mathcal{U}_{2}$, such that
\begin{equation}\label{Eq41}
\left\{
\begin{aligned}
 J_{1}(u_{1}(\cdot), u_{2}(\cdot))=\sup_{v_{1}(\cdot)\in
\mathcal{U}_{1}}J_{1}(v_{1}(\cdot), u_{2}(\cdot)),\\
J_{2}(u_{1}(\cdot), u_{2}(\cdot))=\sup_{v_{2}(\cdot)\in
\mathcal{U}_{2}}J_{2}(u_{1}(\cdot), v_{2}(\cdot)).
\end{aligned} \right.
\end{equation}

\noindent {\bf Solving:} We find the equilibrium point by three steps.

$(i)$ Seek candidate equilibrium points.

Let  $\tilde{q}_{i}(t)$ denote the filtering of $q_{i}(\cdot)$ with respect to $\mathcal{E}_{t}$, i.e. $\tilde{q}_{i}(t)=\mathbb{E}\big(q_{i}(t)\big|\mathcal{E}_{t}\big)$. The similar notations are made for $\tilde{\bar{q}}_{i}(t), \tilde{p}_{i}(t), \tilde{\bar{p}}_{i}(t), \cdots , i=1, 2.$
\begin{align}\label{DG34}
 H_{i}(t, y, z,& Y, Z, v_{1},v_{2},p_{i}, \bar{p}_{i}, q_{i}, \bar{q}_{i})
\triangleq\mbox{\ }\langle
      q_{i}, c_{0}(t)+c_{1}y+c_{2}(t)Y+c_{3}(t)Z\rangle+\langle \bar{q_{i}}, d_{0}(t)\rangle\nonumber\\
&-\langle p_{i}, a_{0}(t)+a_{1}(t)Y+a_{2}(t)Z+F_{i2}(t)v_{1}+a_{4}(t)v_{2}\rangle \mbox{\ }-\langle \bar{p}_{i}, b_{0}(t)\rangle\nonumber\\
& -\frac{1}{2}\Big(\langle e_{i1}(t)y, y\rangle+\langle e_{i2}(t)z, z\rangle+\langle e_{i3}(t)Y, Y\rangle+\langle e_{i4}(t)Z, Z\rangle+\langle e_{i7}(t)v_{i}, v_{i}\rangle \Big).
\end{align}
Applying the maximum principle for nonzero-sum games (Theorem \ref{Theorem2.1}), we confirm that the candidate equilibrium points must satisfy the following form:
\begin{equation}\label{DG30}
    \left\{
    \begin{aligned}
    u_{1}(t)=-e_{17}^{-1}(t)F_{i2}(t)\tilde{p}_{1}(t),\\
    u_{2}(t)=-e_{27}^{-1}(t)a_{4}(t)\tilde{p}_{2}(t),
    \end{aligned}
    \right.
\end{equation}
where $\big(p_{i}(\cdot), \bar{p}_{i}(\cdot), q_{i}(\cdot), \bar{q}_{i}(\cdot)\big)$, for $i=1, 2,$ is the solution of the adjoint FBDSDE:
\begin{equation}\label{DG31}
    \left\{
    \begin{aligned}
        dp_{i}(t)=&\Big(e_{i3}(t)Y(t)+a_{1}p_{i}(t)-c_{2}(t)q_{i}(t)\Big)dt\\
                  &+\Big(e_{i4}(t)Z(t)+a_{2}(t)p_{i}(t)-c_{3}(t)q_{i}(t)\Big)dW(t)-\bar{p}_{i}(t)\hat{d}B(t),\\
        -dq_{i}(t)=&\Big(-e_{i1}(t)y(t)+c_{1}(t)q_{i}(t)\Big)dt-\Big(e_{i2}(t)z(t)\Big)\hat{d}B(t)-\bar{q}_{i}(t)dW(t),\\
        p_{i}(0)=&e_{i6}Y(0)-Mq_{i}(0),\quad q_{i}(T)=-e_{i5}(T)y(T),
    \end{aligned}
    \right.
\end{equation}
and $\big(y(\cdot), z(\cdot), Y(\cdot), Z(\cdot)\big)$ is the solution of the following equation:

\begin{equation}\label{DG32}
\left\{\begin{aligned}
-dY(t)=&\big[a_{0}(t)+a_{1}(t)Y(t)+a_{2}(t)Z(t)-F_{i2}(t)^{2}e_{17}^{-1}(t)p_{1}(t)\\
                      &-a_{4}(t)^{2}e_{27}^{-1}(t)p_{2}(t)\big]dt+b_{0}(t)\hat{d}B(t)-Z(t)dW(t),\\
dy(t)=&\Big[c_{0}(t)+c_{1}(t)y(t)+c_{2}(t)Y(t)+c_{3}(t)Z(t)\Big]dt+d_{0}(t)dW(t)-z(t)\hat{d}B(t),\\
Y(T)=&\ \xi,\;\;\;y(o)=\ MY(0).
\end{aligned}\right.
\end{equation}

$(ii)$ Optimal filtering with the sub-information $\mathcal{E}_{t}=\mathcal{N}\vee \sigma\big\{W(r); 0\leq r\leq t\big\}.$

Equation \eqref{DG31} together with \eqref{DG32} constitutes a triple dimensional FBDSDE. In order to find the explicit expression of the candidate equilibrium point, we need to compute the optimal filters
$\tilde{p}_{1}(\cdot)$ and $\tilde{p}_{2}(\cdot)$ of $p_{1}(\cdot)$ and $p_{1}(\cdot)$.  Applying the filtering result derived by Xiong (\cite{Xiong2008}, Lemma 5.4) to \eqref{DG31} and \eqref{DG32} under the available sub-information $\mathcal{E}_{t}=\mathcal{N}\vee \sigma\big\{W(r); 0\leq r\leq t\big\}$, we conclude that $\tilde{p}_{1}(\cdot)$ and $\tilde{p}_{2}(\cdot)$ satisfy the following triple dimensional FBSDE:
\begin{equation}\label{DG33}
    \left\{
      \begin{aligned}
        -d\left(\begin{array}{c}
           \tilde{Y}(t)\\
           \tilde{q}_{1}(t)\\
           \tilde{q}_{2}(t)\\
           \end{array}
                       \right)=&\left\{\left(\begin{array}{ccc}
                                     a_{1}(t) & 0 & 0 \\
                                     0        &  c_{1}(t)                      & 0\\
                                     0        &  0                              & c_{1}(t)\\
                                    \end{array}
                                                \right)\left(\begin{array}{c}
                                                          \tilde{Y}(t)\\
                                                         \tilde{q}_{1}(t)\\
                                                        \tilde{q}_{2}(t)\\
                                                         \end{array}
                                                                         \right)\right.\\
                               &\hspace{4mm}+\left(\begin{array}{ccc}
                                     0 & -F_{i2}(t)^{2}e_{17}^{-1}(t) & -a_{4}(t)^{2}e_{28}^{-1}(t) \\
                                     -e_{11}(t) &    0                    & 0\\
                                     -e_{21}(t) &  0                     & 0\\
                                    \end{array}
                                                \right)\left(\begin{array}{c}
                                                          \tilde{y}(t)\\
                                                         \tilde{p}_{1}(t)\\
                                                        \tilde{p}_{2}(t)\\
                                                         \end{array}
                                                                         \right)\\
                            &\hspace{4mm}+\left(\begin{array}{ccc}
                                     a_{2}(t) & 0 & 0 \\
                                     0 &    0                    & 0\\
                                    0 &  0                     & 0\\
                                    \end{array}
                                                \right)\left(\begin{array}{c}
                                                          \tilde{Z}(t)\\
                                                         \tilde{\bar{q}}_{1}(t)\\
                                                        \tilde{\bar{q}}_{2}(t)\\
                                                         \end{array}
                                                                         \right)+\left.\left(
                                     \begin{array}{c}
                                       a_{0}(t)\\
                                       0\\
                                       0
                                     \end{array}
                              \right)\right\}dt-\left(
                                                    \begin{array}{c}
                                                      \tilde{Z}(t)\\
                                                      \tilde{\bar{q}}_{1}(t)\\
                                                      \tilde{\bar{q}}_{2}(t)
                                                      \end{array}
                                            \right)dW(t)\\
       d\left(\begin{array}{c}
           \tilde{y}(t)\\
           \tilde{p}_{1}(t)\\
           \tilde{p}_{2}(t)\\
           \end{array}
                       \right)=&\left\{\left(\begin{array}{ccc}
                                     c_{2}(t) & 0 & 0 \\
                                     e_{13}(t) &  -c_{2}(t)                       & 0\\
                                     e_{23}(t) &  0                     & -c_{2}(t)\\
                                    \end{array}
                                                \right)\left(\begin{array}{c}
                                                          \tilde{Y}(t)\\
                                                         \tilde{q}_{1}(t)\\
                                                        \tilde{q}_{2}(t)\\
                                                         \end{array}
                                                                         \right)\right.\\
                               &\hspace{4mm}+\left(\begin{array}{ccc}
                                     c_{1}(t) & 0 & 0 \\
                                     0        &  a_{1}(t)                      & 0\\
                                     0        &  0                              & a_{1}(t)\\
                                    \end{array}
                                                \right)\left(\begin{array}{c}
                                                          \tilde{y}(t)\\
                                                         \tilde{p}_{1}(t)\\
                                                        \tilde{p}_{2}(t)\\
                                                         \end{array}
                                                                         \right)\\
                             &\hspace{4mm}+\left(\begin{array}{ccc}
                                     c_{3}(t) & 0 & 0 \\
                                     0        &  0 & 0\\
                                     0        &  0  & 0\\
                                    \end{array}
                                                \right)\left(\begin{array}{c}
                                                          \tilde{Z}(t)\\
                                                         \tilde{\bar{q}}_{1}(t)\\
                                                        \tilde{\bar{q}}_{2}(t)\\
                                                         \end{array}
                                                                         \right)+\left.\left(
                                     \begin{array}{c}
                                       c_{0}(t)\\
                                       0\\
                                       0
                                     \end{array}
                               \right)\right\}dt\\
                             &+\left\{\left(\begin{array}{ccc}
                                     0 & 0 & 0 \\
                                    0 &  -c_{3}(t)                       & 0\\
                                    0 &  0                     & -c_{3}(t)\\
                                    \end{array}
                                                \right)\left(\begin{array}{c}
                                                          \tilde{Y}(t)\\
                                                         \tilde{q}_{1}(t)\\
                                                        \tilde{q}_{2}(t)\\
                                                         \end{array}
                                                                         \right)
                                    +\left(\begin{array}{ccc}
                                     0 & 0 & 0 \\
                                    e_{14}(t) &  0                       & 0\\
                                    e_{24}(t)  &  0                     & 0\\
                                    \end{array}
                                                \right)\left(\begin{array}{c}
                                                          \tilde{Z}(t)\\
                                                         \tilde{\bar{q}}_{1}(t)\\
                                                        \tilde{\bar{q}}_{2}(t)\\
                                                         \end{array}
                                                                         \right)\right.\\
                            &\hspace{4mm}+\left(\begin{array}{ccc}
                                     0 & 0 & 0 \\
                                     0        &  a_{2}(t)                      & 0\\
                                     0        &  0                              & a_{2}(t)\\
                                    \end{array}
                                                \right)\left(\begin{array}{c}
                                                          \tilde{y}(t)\\
                                                         \tilde{p}_{1}(t)\\
                                                        \tilde{p}_{2}(t)\\
                                                         \end{array}
                                                                         \right)+\left.\left(
                                                                                 \begin{array}{c}
                                                                                    d_{0}(t)\\
                                                                                      0\\
                                                                                       0
                                                                                  \end{array}
                                                                                   \right)\right\}dW(t),\\
    \left(
         \begin{array}{c}
         \tilde{Y}(T)\\
         \tilde{q}_{1}(T)\\
         \tilde{q}_{2}(T)\\
          \end{array}
          \right)=&\left(
                          \begin{array}{c}
                           \mathbb{E}[\xi|\mathcal{E}_{T}]\\
                           -e_{15}(T)\tilde{y}(T)\\
                           -e_{25}(T)\tilde{y}(T)\\
                           \end{array}
                           \right),
          \left(
         \begin{array}{c}
         \tilde{y}(0)\\
         \tilde{p}_{1}(0)\\
         \tilde{p}_{2}(0)\\
          \end{array}
          \right)=\left(\begin{array}{ccc}
                                     M & 0 & 0 \\
                                    e_{16}(t) &  -M                       & 0\\
                                    e_{26}(t)&  0                     & -M\\
                                    \end{array}
                                                \right)\left(\begin{array}{c}
                                                          \tilde{Y}(0)\\
                                                         \tilde{q}_{1}(0)\\
                                                        \tilde{q}_{2}(0)\\
                                                         \end{array}
                                                                         \right).
\end{aligned}
    \right.
\end{equation}
Just like Huang et al.\cite{HWX2009}, we call \eqref{DG33} a forward-backward stochastic differential filtering equation, which is distinguished from the classical filtering literature
(see e.g. Liptser and Shiryaev\cite{LS1977}, Xiong\cite{Xiong2008}).  Now, we obtain an explicit candidate equilibrium point for the foregoing LQ nonzero-sum differential game.

 $(iii)$ Verify that $\big(u_{1}(\cdot), u_{2}(\cdot)\big)$ denoted by \eqref{DG30} is indeed an equilibrium point.

We can check that the system \eqref{DG19} and performance criterion \eqref{DG20} satisfy the assumptions (H1) and (H2),  the Hamiltonian $H_{i}$ $(i=1, 2)$ denoted by \eqref{DG34} satisfies the conditions (\ref{Eq73}--\ref{Eq60'}). Then, from Theorem \ref{Thm3.2}, we conclude that $\big(u_{1}(\cdot), u_{2}(\cdot)\big)$ denoted by \eqref{DG30} is indeed an equilibrium point.

\section{Conclusion}

We are concerned with a new type of stochastic differential game problems of FBDSDEs. There are two distinguishing features: One is that  game systems are initial coupled; The other one is that differential games is under partial information.  We established a  maximum principle and a verification theorem, also called a necessary condition and a sufficient condition, for an equilibrium point of partial information nonzero-sum stochastic differential games. Zero-sum games can be considered as a particular case of nonzero-sum games, so we also gave the corresponding conditions for a saddle point of zero-sum stochastic differential games. Finally, we worked out an LQ example and gave the explicit expression of an equilibrium point of nonzero-sum differential games.

It is worth pointing out that game system of FBDSDE covers many cases as its particular case. If we drop its the terms of backward It\^{o}'s integral or forward equation or both them , FBDSDE can be reduced to FBSDE or BDSDE or BSDE. Moreover, if we suppose that $\mathcal{E}_{t}=\mathcal{F}_{t}$ for all $t\in [0, T],$  all the results are reduced to the case of full information. In addition, stochastic control problems can be regarded as zero-sum stochastic differential games with only one player. Then, our results are a partial extension to Xiao and Wang\cite{XWJAMC} for optimal control of FBSDEs with partial information, Han et al.\cite{HPW2010} for optimal control of BDSDEs with full information, Huang et al.\cite{HWX2009} for optimal control of BSDEs with partial information, Wang and Yu\cite{WY2010} and Yu and Ji\cite{YJ2008} for differential games of BSDEs with full information, and Wang and Yu\cite{WY2011} for differential games of BSDEs with partial information. In our game systems of FBDSDEs with partial information, the forward equations are coupled with the backward equations at initial time, not terminal time. So they do not cover each other between our results and those of terminal coupled forward-backward stochastic systems with full or partial information derived by Buckdahn and Li\cite{BL2008}, Hamad\`{e}ne\cite{Ha1999}, Hui and Xiao\cite{HX2011}, Meng\cite{Meng2009,Meng2010}, {\O}ksendal and Sulem\cite{OS2010}, Peng and Wu\cite{PW1999}, Shi and Wu\cite{SW2006,SW-appear},Wang and Wu\cite{WW2008,WW2009}, Wu\cite{Wu2005,Wu2010}, Xiao and Wang\cite{XWSAA}, Zhang and Shi\cite{ZS2010}, Zhu et al.\cite{ZSG2009}.

Finally, since there are many partial information optimization and game
problems in finance and economics, we hope that the results have
applications in these areas.

\end{document}